\documentclass[letterpaper, 11pt, reqno]{amsart}

\usepackage{setspace}\usepackage{mathrsfs}
\usepackage[utf8]{inputenc}
\usepackage{graphicx}\usepackage{constants}
\usepackage{amsmath}
\usepackage{amssymb,amsthm}
\usepackage{cases}
\usepackage{esint}
\usepackage[a4paper,margin=3cm]{geometry}
\usepackage[citecolor=blue, linkcolor=blue,colorlinks]{hyperref}

\graphicspath{ {images/} }
\newtheorem{theorem}{Theorem}
\newtheorem{lemma}{Lemma}

\newtheorem{prop}{Proposition}
\newtheorem{de}{Definition}

\theoremstyle{definition}
\newtheorem{rk}{Remark}

\DeclareMathOperator*{\supp}{supp}

\def\R{\mathbb{R}}

\def\half{{\textstyle\frac{1}{2}}}

\def\a{\alpha}
\def\e{\epsilon} 

\def\p{\partial}
\def\la{\lambda}
\def\k{\kappa}\font\smathbold=msbm7\font\ssmathbold=msbm7 at 5.5pt
\def\sR{{\hbox{\smathbold\char82}}}\def\ssR{{\hbox{\ssmathbold\char82}}}
\def\ds{\displaystyle}\def\Br{B_{\sR^n}}\def\xb{{\overline x}}\def\yb{{\overline y}}
\def\ni{\noindent}

\def\ba{\begin{aligned}}
\def\ea{\end{aligned}}

\def\wt{\widetilde}\def\L{\Lambda}
\def\b{\beta}\def\dist{{\rm dist}}

 \def\cut{{\text {cut}}}
 
\def\be{\begin{equation}}\def\ee{\end{equation}}\def\nt{\notag}\def\bc{\begin{cases}}\def\ec{\end{cases}}
\def\inj{{\text {inj}}}\def\N{{\mathbb N}}
\def\qed{\rightline{\setlength{\fboxsep}{0pt}\setlength{\fboxrule}{0.2pt}\fbox{\rule[0pt]{0pt}{1.3ex}\rule[0pt]{1.3ex}{0pt}}}}\def\cut{{\text {Cut}}}

\def\nna{{\frac{n}{n-\a}}}\def\na{\frac{n}{\a}}
\def\ric{{\rm{Ric } }}

\def\llceil{\left\lceil}\def\rrceil{\right\rceil}
\newcommand\fH[1]{\sbox0{#1}\dimen0=\ht0 \advance\dimen0 -1ex
  \sbox2{\'{}}\sbox2{\raise\dimen0\box2}%
  {\ooalign{\hidewidth\kern.1em\copy2\kern-.5\wd2\box2\hidewidth\cr\box0\crcr}}}\def\rb{ \overline r}
\def\S{\mathbb{S}}
 
\font\twelvemi=cmmi12 at 12pt\font\elevenmi=cmmi11 at 9 pt\font\fivemi=cmmi5 at 6 pt\font\twelvemimi=cmmi12 at 12pt
\renewcommand{\chi}{\raisebox{.13\baselineskip}{\hbox{\twelvemi\char31}}}
\newcommand{\gam}{{\raisebox{.08\baselineskip}{\hbox{\twelvemi\char13}}}}
\newcommand{\sgam}{{\raisebox{.08\baselineskip}{\hbox{\elevenmi\char13}}}}
\newcommand{\sm}{{\hbox{\twelvemimi\char27}}_{\!\hbox{\fivemi\char 77}}}

\newcommand{\sx}{{\raisebox{.0\baselineskip}{\hbox{\twelvemimi\char27}}_{\!x}}}
\def\an{\frac{\a}{n}}
\renewcommand{\gamma}{\gam}
\def\lan{\big\langle}\def\ran{\big\rangle}\def\vol{{\text {Vol}}}\def\quarter{{\textstyle \frac1 4}}
 \def\Ar{{A_{r,R}}}\def\smi{\sm^{-1}}\def\bN{{\bf N}}\def\mut{\tilde\mu}
\def\Gt{\widetilde G} \def\ct{\widetilde c}\def\Grn{{G_\a^{\sR^{\hskip-.007in n}}}}\def\Drn{\Delta_{\sR^{\hskip-.007in n}}}
\vskip 0.1in\def\grn_#1{{G_#1^{\sR^{\hskip-.007in n}}}}\def\sgrn_#1{{G_#1^{\ssR^{\hskip-.007in n}}}}
\def\per{{\text {Per}}}\def\Lt{\wt\L}
\def\tx{\theta_x}

\title[Sharp inequalities on manifolds with Euclidean volume growth] {Sharp estimates and inequalities on Riemannian manifolds with Euclidean volume growth}
\author[L. Fontana, C. Morpurgo, L. Qin] {Luigi Fontana, Carlo Morpurgo,  Liuyu Qin }
\thanks{The third author was supported by the National Natural Science Foundation of China (12201197) }
\begin{document}
\begin{abstract} We obtain sharp estimates for heat kernels and Green's functions on  complete noncompact Riemannian manifolds with Euclidean volume growth and nonnegative Ricci curvature. 
 We will then apply these estimates to obtain sharp Moser-Trudinger inequalities on such manifolds. 

\end{abstract}
\numberwithin{equation}{section}

\maketitle

 \section{Introduction}

A complete, noncompact, $n-$dimensional  Riemannian manifold with nonnegative Ricci curvature is said to have \emph{Euclidean volume growth}, or EVG,  if 
\be \frac{\vol\big(B(x,r)\big)}{|B_{\sR^n}(r)|}\searrow \sm>0,\qquad {\text {as }}\;r\to \infty\label{evg0}\ee
where $\vol(B(x,r)$ denotes the volume of the  geodesic ball centered at $x$ and radius $r$,  and $|B_{\sR^n}(r)|$ the volume of the Euclidean ball of radius $r$.  The volume ratio in \eqref{evg0} is decreasing in $r$ by the Bishop comparison theorem, which also guarantees that $\vol\big(B(x,r)\big)\le{|B_{\sR^n}(r)|}$, and therefore $0<\sm\le1.$ It's easy to see that the limit in \eqref{evg0} is independent of $x$.

Manifolds with EVG have been extensively studied since the late 80s, see for example \cite{bkn}, \cite{cm1}, \cite{ct}, \cite{s}, \cite{xia1}, \cite{xia2}.

Very recently, several authors found  versions of the classical isoperimetric and Sobolev inequalities on manifolds with EVG, whose sharp constants feature explicitly the  constant $\sm$, known currently as \emph{Asymptotic Volume Ratio}, or AVR.
In particular, the sharp isoperimetric inequality holds in the form 
\be \per(\Omega)\ge \sm^{\frac1 n} nB_n^{\frac1 n}\vol(\Omega)^{\frac{n-1}n}\label{isop}\ee
where $B_n=|B_{\sR^n}(1)|$ and where $\per(\Omega)$ denoted the perimeter of a bounded open set $\Omega$ with smooth boundary. This result, along with  the characterization of equality,  was proved by Brendle \cite{br}, while other authors proved the result for more general $\Omega$, and on more general metric measure spaces, with different techniques (see \cite{bk}, \cite{apps}, \cite{afm}).

In \cite{bk}, the sharp Sobolev inequality on manifolds with EVG was derived in the form
\be \|u\|_{p^*}\le \sm^{-\frac 1 n} C(n,p)\|\nabla u\|_p,\qquad 1<p<n,\;\; p^*=\frac{np}{n-p}\label{sobolev}\ee
where $C(n,p)$ denotes the best constant in the Euclidean version of the same inequality \cite{au}, \cite{t2}.

In this paper we will focus on the borderline case of the Sobolev embedding theorem, which in $\R^n$ concerns the Sobolev spaces $W^{\a,\frac n\a}$. In particular we will obtain sharp  Moser-Trudinger inequalities on manifolds with EVG. In order to state our main results we will need a few more definitions (for some more details see the next section).

\smallskip
Let us define  
\be D^\a=\bc (-\Delta)^{\frac{\a}{2}}  &\text{  for $\alpha$ even}\cr\nabla(-\Delta)^{\frac{\a-1}{2}} & \text{ for $\alpha$ odd.} \label{Da}
\ec\ee
where $\Delta$ is the Laplace-Beltrami operator on $M$ and $\nabla$ is the gradient. 

The Bessel potential space is defined for $\a>0,p\ge1$  as 
\be H^{\a,p}(M)=\{u\in L^p(M): \; u=(I-\Delta)^{-\frac\a2}f,\; f\in L^p(M)\}\label {bessel1}\ee
endowed with the norm $\|(I-\Delta)^{\frac\a2}u\|_p$. These spaces were studied, on general noncompact manifolds, by Strichartz \cite{st}, who also showed  that $C_c^\infty(M)$ is dense in  $H^{\a,p}(M)$ and that $\|u\|_p+\|(-\Delta)^{\frac\a2}u\|_p$ is an equivalent norm if $p>1$. By work of Auscher et al \cite{acdh}, it turns out that in our setting  the Riesz transform is continuous in $L^p$ and this implies that $\|u\|_p+\|\nabla(-\Delta)^{\frac{\a-1}2}u\|_p$
is an equivalent norm in $H^{\a,p}$, for $\a\ge1.$ See Section 2 for some comments about the relations between $H^{\a,p}$ and the classical Sobolev spaces $W^{\a,p}$ and $W_0^{\a,p}.$

If  $\Drn$ is the standard Laplacian on $\R^n$ then the fundamental solution of  $(-\Drn)^{\frac{\a}{2}}$ will be denoted by 
\be \Grn(\xi,\eta)=N_\a(\xi-\eta),\qquad N_\a(\xi)=c_\a|\xi|^{\a-n}\label{N}\ee
i.e. the classical Riesz potential, and the explicit value of $c_\a$ is given in \eqref{GE}.

Throughout this paper the volume of the unit ball in $\R^n$ will also be denoted by $B_n$:
\be B_n=|B_{\sR^n}(1)|=\frac{\omega_{n-1}}{n}\ee
where $\omega_{n-1}=2\pi^{n/2}/\Gamma(n/2)$ is the surface area of the $(n-1)$-dimensional sphere $\S^{n-1}=\partial \Br(1)$.

The \emph{Adams constant} is defined as
 \be \gam_{n,\a}=
\begin{cases}
\dfrac{c_\a^{-\frac{n}{n-\a}}}{B_n}\ \ \ &\text{if} \ \a \ \text{even}\\[1em]
\dfrac{\ct_{\a}^{\hskip .21em-\frac{n}{n-\a}}}{B_n}\ \ &\text{if} \  \a\ \ \text{odd},
\end{cases}
\label{gamma}
\ee
with $\ct_\a$ defined in \eqref{GTE}, (related to $c_\a$ as in \eqref{GTE1}),
and which gives the sharp exponential constant in the classical Moser-Trudinger inequality on bounded domains of $\R^n$, and in  many other related inequalities.
 
 Finally, the \emph{regularized exponential function} is defined as
\be \exp_m(t)=e^t-\sum_{k=0}^{m}\frac{t^k}{k!},\qquad t\ge0,\; m=0,1,2...\ee
 
  \begin{theorem}\label{main1} Let $(M,g)$ be a complete, connected, noncompact Riemannian manifold with  nonnegative Ricci curvature and Euclidean volume growth.

For  $\a$ even and $\a<n$, or for $\alpha$ odd and $\a\le n/2$,   there exists $C>0$ such that for all $u\in H^{\a,\na}(M)$
with 
\be \|D^\a u\|_{n/\a}\le 1\label {norm1}\ee we have 
\be  \int_M \frac{\exp_{\llceil\frac{n-\a}{\a}-1\rrceil}{\left(\sm^{\frac{\a}{n-\a}}\gam_{n,\a}|u|^{\frac{n}{n-\a}}\right)}}{1+|u|^{\nna}}d\mu\le C\|u\|_{{ n/\a}}^{{n/\a}}. \label{MS}\ee

If $\alpha=1$, $\alpha=2$,   $4\le \a<\frac n2$ and $\alpha$ even, $3\le \alpha<\frac n 2+1$ and  $\alpha$ odd,  then the exponential constant in  \eqref{MS} is sharp, i.e. it cannot be replaced by a larger value. For the same range of $\alpha$, the  power $\nna$ in the denominator in  \eqref{MS} is sharp, i.e. it cannot be replaced by a smaller number.  
  \end{theorem}
  \smallskip  Theorem \ref{main1} includes $\R^n$ with the standard metric, in which case $\sm=1$,  where the inequalities hold for all $\alpha<n$ (see \cite{imn}, \cite{ms1}, \cite{ms2}, \cite{ms3}, \cite{ltz}, \cite{q}). 
  
  On noncompact, simply connected Riemannian manifolds with nonpositive sectional curvature $K$, i.e. Hadamard manifolds, \eqref{MS} holds (for all $\a<n$) with sharp exponential constant $\gamma_{n,\a}$, provided $-b^2\le K\le -a^2<0$ \cite{mq1}, \cite{mq2} (see also \cite{mq1} for a sharp version of \eqref{MS} on the Heisenberg group).
 
 We emphasize that all of the results known so far, on the above  Riemannian manifolds,  hold with sharp constant $\gamma_{n,\a}$, whereas within manifolds with EVG the sharp constant is \emph{lower}, by a factor of $\sm^{\frac{\a}{n-\a}}$. The main point is that inequalities such as \eqref{MS} under the  norm condition \eqref{norm1} are very sensitive to the geometry of the ambient space at infinity, in particular volume growth,  something that  had already been observed in \cite{mq1} (see for example the exponential constant in the Adams inequality \eqref{1a} below).
 
 \smallskip
 
 In the next theorem, we will show  that on manifolds with EVG, and  with enough  bounded geometry,  imposing a  norm condition stronger than \eqref{norm1}  yields  sharp Moser-Trudinger inequalities with the usual, classical  constant $\gamma_{n,\a}$.
 
 \smallskip 
 
We will say that  $M$ has \emph{$k^{\text {th}}$ order bounded geometry} if there is $B>0$ such that 
\be \inj(M)>0,  
\qquad \quad |\nabla^i R|\le B, \qquad  i=0,1,...,k.\label{M4} \ee
where $\inj(M)$ is the injectivity radius of $M$, $R$ the curvature tensor, and where $\nabla^i$ is the $i^{th}$ covariant derivative.

 \begin{theorem}\label{main2} Let $(M,g)$ be a complete, connected, noncompact Riemannian manifold with  nonnegative Ricci curvature and Euclidean volume growth.  
 For any $\alpha$ integer with $1\le\alpha<n$,  let $M$ have $1^{st}$ order bounded geometry for $\alpha$ even, and  $2^{nd}$ order bounded geometry for $\a$ odd.  Then for any $\kappa>0$ there is $C$ such that for  all $u\in H^{\a,\na}(M)$
with 
\be\kappa \|u\|_{n/\a}^{n/\a}+ \|D^\a u\|_{n/\a}^{n/\a}\le 1\label {norm2}\ee we have 
\be  \int_M \frac{\exp_{\llceil\frac{n-\a}{\a}-1\rrceil}{\left(\gam_{n,\a}|u|^{\frac{n}{n-\a}}\right)}}{1+|u|^{\nna}}d\mu\le C\|u\|_{{ n/\a}}^{{n/\a}}, \label{MSR}\ee
and
\be  \int_M \exp_{\llceil\frac{n-\a}{\a}-1\rrceil}{\left(\gam_{n,\a}|u|^{\frac{n}{n-\a}}\right)}d\mu\le C. \label{MT}\ee 
 The exponential constant in \eqref{MSR}  and \eqref{MT} is sharp, and the power $\nna$ in the denominator in \eqref{MSR} is sharp.
  
 \end{theorem}

 Inequality  \eqref{MT} on $\R^n$, was proved in \cite{ruf}, \cite{lr}, \cite{ll}, \cite{fm2}. It is easy to check that  on $\R^n$   inequality \eqref{MT} under \eqref{norm2} (with $\kappa>0$) follows from the same result with $\kappa=1$,  by a  dilation argument. It is worth mentioning, at this point, that  \eqref{MT} is \emph{false} under \eqref{norm1}, that is under the weaker condition $\|D^\alpha u\|_{n/\a}\le 1.$ In fact, \eqref{MT} is false even under the full Sobolev norm conditions
 \be \max\big\{\|u\|_{n/\a},\|D^\a u\|_{n/\a}\big\}\le1 \qquad  \big(\|u\|_{n/\a}^{qn/\a}+\|D^\a u\|_{n/\a}^{qn/\a}\big)^{\frac\a{qn}}\le 1,\quad q>1. \ee
 Under such  conditions one can only prove \eqref{MT} with $\gam_{n,\a}$ replaced by any $\gam<\gam_{n,\a}$. In particular, for any $\gamma<\gamma_{n,\a}$, there is $C$ (depending on $\gamma$) such that under \eqref{norm1} there holds
 \be  \int_{\R^n} \exp_{\llceil\frac{n-\a}{\a}-1\rrceil}{\left(\gam|u|^{\frac{n}{n-\a}}\right)}d\mu\le C \|u\|_{n/\a}^{n/\a}.\label{AT}\ee
This  inequality is often called of Adachi-Tanaka type,  since it first appeared when $\alpha=1$ in a paper by Adachi and Tanaka [AT];  later it was extended to all integer values of $\alpha$ in \cite{fm2},
 as a relatively easy consequence of \eqref{MT} under  \eqref{norm2}. In \cite{q} Corollary 1, Qin derived \eqref{AT} as an easy consequence of \eqref{MSR} under either \eqref{norm1} or \eqref{norm2}  when $u$ is a  general Riesz-like potential. For a short discussion and history of  such types of weaker,  \emph{subcritical} Moser-Trudinger inequalities, which will not be treated in this paper,  see \cite{fm2} and references therein (in particular, see \cite{fm2} Corollary 4). 
 
 Still on $\R^n$,  when $\gamma=\gamma_{n,\a}$ inequality  \eqref{AT}  fails  under \eqref{norm1}, however it becomes true, and sharp, if the integrand is divided by $1+|u|^\nna$, which gives  \eqref{MSR} (and which is the same as the already mentioned \eqref{MS}, since $\sm=1$ on $\R^n$). This was first discovered in \cite{imn}, when $n=2,\alpha=1$,  followed by the aforementioned papers  \cite{ms1}, \cite{ms2}, \cite{ms3}, \cite{ltz}, \cite{q}, which treated other values of $n$ and $\alpha$.  For a unified discussion of all the known most important Moser-Trudinger inequalities on domains of $\R^n$, see the introduction of \cite{mq1}.

  It has been remarked before that on $\R^n$  inequality \eqref{MS} (which is the same as \eqref{MSR}) under \eqref{norm1} is stronger than  \eqref{MT} under   \eqref{norm2}, in the sense that it directly implies it via H\"older's inequality and \eqref{AT}  (see \cite{ms1}, \cite{ms2}, \cite{ms3}, and also  \cite{q}, Section 8).  With the same arguments one sees, on general manifolds,  that  if \eqref{MSR} holds with exponential constant $\gamma^*$ under \eqref{norm1},  then  \eqref{AT} also holds for $\gamma<\gamma^*$, and \eqref{MT} holds,  under \eqref{norm2}, with the \emph{same} exponential constant~$\gamma^*$.  
With small modifications of the arguments   in \cite{q} Section 8, the same is true if one assumes \eqref{MSR} with exponential constant $\gamma^*$ under \eqref{norm2}, rather than \eqref{norm1}; in particular, in the context of Theorem 2, estimate \eqref{MSR} implies  estimate \eqref{MT} (however we do not know if the converse implication is true). 
  One consequence of all this is that on manifolds with EVG and enough bounded geometry, inequality \eqref{MS} under \eqref{norm1} is \emph{not} stronger than \eqref{MT} under \eqref{norm2}, unless $\sm=1$, i.e. $M=\R^n$, with its standard metric. On such manifolds, therefore, the sharp estimates contained in Theorems 1 and 2 are in a sense independent of one another.

On Hadamard manifolds,   \eqref{MT} holds when $\alpha=1$ under \eqref{norm2},  provided $-b^2\le K\le0$ \cite{fmq} (with the exception $n=2,3,4$, in which cases it holds for $K\le0$ \cite{kr}). When $-b^2\le K\le -a^2<0$, estimate \eqref{MT} holds under the weaker condition \eqref{norm1} (\cite{bs}). To the best of our knowledge the results in Theorem \ref{main2}, even when $\alpha=1$,  are not known on other noncompact manifolds with variable  Ricci curvature which is  only bounded below (or even nonnegative),  even assuming positive injectivity radius (see \cite{kr}, \cite{y} for related results).

 The results in Theorem \ref{main2}, and their  proofs, are more local in nature, and somewhat insensitive to the volume growth of large balls. In particular, in \cite{fmq} the authors introduced the concept of ``Local Moser-Trudinger inequality'', for the case $\alpha=1$, and showed that if an inequality like \eqref{MT} holds for functions whose support have volume less than a given number, then the inequality holds globally. In fact, it is  possible to prove \eqref{MT} above using this result, for the specific case $\alpha=1$ (see Remark \ref{LMT}).

 There are essentially two  methods of proof of basic Moser-Trudinger type inequalities. In the first approach, used mainly in $\R^n$,  one tries to replace the original function $u$  with its symmetric decreasing rearrangement $u^\#$, thereby reducing the problem to a sharp one-dimensional estimate.  On $\R^n$, and for $\alpha=1$, this is done via the P\'olya-Szeg\H o inequality $\|\nabla u^\#\|_p\le \|\nabla u\|_p$ (see for example \cite{mo}, \cite{ruf}, \cite{lr}, \cite{imn}, \cite{ms1}, \cite{ms2}). For $\alpha\ge2$ one can use Talenti's comparison theorem, to establish more general inequalities for Sobolev  functions satisfying Navier boundary conditions see \cite{ms3}, \cite{t}.

 The second approach is to replace $u$ by a suitable potential, typically after  an estimate of type
 \be |u(x)|\le \int_M |K_\a(x,y)|\,|D^\a u(y)| d\mu(y)\label{potential}\ee
 where $K_\a$ is essentially a fundamental solution of $D^\a$. Indeed, one  typically takes  $K_2(x,y)=G(x,y)$, the Green function of $\-\Delta$, and $K_1(x,y)=\nabla_y G_(x,y)$. Through sharp asymptotic estimates of $K_\a$ along the diagonal, and at infinity for some noncompact manifolds, one is generally able to follow the technique originally developed by Adams on $\R^n$ in \cite{a1}, and to derive a sharp exponential inequality where $u$ is replaced by the integral operator with kernel $K_\a$. This method has been successfully used in spaces where the symmetrization results mentioned above  are not (yet)  available, such as compact manifolds, CR manifolds, Hadamard manifolds etc. (see for example \cite{a1}, \cite{f},  \cite{fm1}, \cite{fm2}, \cite{q}, \cite{mq1}).
 
 Both methods have shortcomings on general noncompact manifolds satisfying $\inj(M)>~0$  and $\ric\ge~\lambda $. 
 The P\'olya-Szeg\H o inequality holds in the form 
 \be I(\Omega) \|\nabla u^\#\|_{L^p(\Omega^\#)}\le  \|\nabla u\|_{L^p(\Omega)},\qquad u\in W_0^{1,p}(\Omega)\label {PS}\ee
where $\Omega$ is an open bounded set in $M$, $\Omega^\#$ a ball  in $\R^n$ with volume $\vol(\Omega)$,  and where $I(\Omega)$ is the normalized isoperimetric ratio
\be I(\Omega)=\inf_U \frac{\per(U)}{nB_n^{\frac1 n}\vol(U)^{\frac{n-1}n}}\in[0,1]
\ee 
 the infimum being taken among smooth open subsets of $\Omega.$ Assuming $I(\Omega)>0$,  this result will allow easily to prove a Moser-Trudinger inequality for functions in $W_0^{1,n}(\Omega)$, with exponential constant $I(\Omega)^{\frac n{n-1}}\gamma_{n,1}$ (see for ex. \cite{kr}). However, proving the sharpness of such constant is not possible without further information on the sharpness of the isoperimetric inequality, and ultimately on the geometry of $M$.
 
Secondly, establishing uniform asymptotic estimates for $K_\a$ is not always  an easy task. For $\alpha$ even, if $K_\a(x,y)=G_\a(x,y)$ is a fundamental solution of $(-\Delta)^{\frac\a2}$, the ``Adams' method" requires at the very least a local asymptotic estimate of type 
 \be |G_\a(x,y)- c_\a d(x,y)^{\a-n}|\le C d(x,y)^{a-n+\e}  ,\qquad d(x,y)\le 1\label{e1}\ee
 with  $C$ \emph{independent of $x,y$} (it is not difficult to show \eqref{e1} for fixed $x$ and with $C$ depending on $x$). To prove such result one needs to impose further assumptions on the curvature.
 
Moreover,  it is necessary to have some control on the decay of $G_\a(x,y)$ at infinity, and for the validity of an inequality  such as \eqref{MS} it is also important to establish more precise  asymptotic bounds of $G(x,y)$ for large distances.

On manifolds with EVG we are able to prove \eqref{MS} with a technique which is a sort of hybrid between the two methods discussed above. One one hand, there is a sharp  P\'olya-Szeg\H o inequality
\be \|\nabla u^\#\|_{L^p(\Omega^\#)}\le \sm^{-\frac1n}\|\nabla u\|_{L^p(\Omega)},\qquad u\in W_0^{1,p}(\Omega)\label{PS1}\ee
(see \cite{bk}, \cite{apps}), which readily implies \eqref{MS}  for $\alpha=1$ (but not the sharpness part!); see \eqref{PSZ} and related comments. On the other hand, to settle the case $\alpha>1$, in Section \ref{Tal} we will  obtain sharp Talenti/Weinberger type symmetrization estimates on the fundamental solutions of $D^\a$, by using recent results by Chen-Li \cite{cl}.  For $\alpha$ even we will obtain the estimate 
\be \big(G_\a(x,\cdot)\big)^*(t)\le \sm^{-\frac\a n}N_\a^*(t),\qquad t>0\ee
for all $x\in M$, where $f^*$ denotes the  decreasing rearrangement of a function $f$, and where $N_\a(\xi)$ is the Riesz kernel on $\R^n$, defined in \eqref{N}. We will  also  obtain a sharp comparison result for 
\be\mathop\int\limits_{s_1<G_\a(x,y)\le s_2} |\nabla_y G_\a(x,y)|^q \;d\mu(y) \ee
in terms of $N_\a$, and for $q\le 2$ (see Theorem \ref{talenti}). These results will allow us to get one of the key estimates needed to obtain \eqref{MS}, namely
\be \mathop\int\limits_{r_1<d(x,y)\le r_2} |K_\a(x,y)|^{\frac n {n-\a}}\; d\mu(y)\le \sm^{-\frac\a{n-\a}}\gamma_{n,\a}^{-1}\log\frac{V_x(r_1)}{V_x(r_2)}+B,\qquad 0<r_1<r_2\ee
where $K_\a=G_\a$ for $\alpha$ even and $K_\a=\nabla_y G_{\a+1}$ for $\alpha$ odd, in which case we require $\nna\le 2$, i.e. $\a\le n/2$. The inequality in \eqref{MS} will then follow from the representation formula \eqref{potential}, combined with an Adams inequality on metric measure spaces, Theorem~\ref{m1}, recently obtained  in \cite{mq1}.  At this point  we do not know  how to handle the case $\alpha$ odd and $n/2<\a<n$.

The proof of \eqref{MSR} and \eqref{MT}  will instead be based on the uniform local estimate \eqref{e1}, and a similar one for $|\nabla_y G_{\a+1}|$, and not on the sharp behavior of the fundamental solution at infinity. These estimates will be obtained in Section \ref{greensmall} as a consequence of local, uniform, small times  estimates for the heat kernel, which we will derive in Section~\ref{heatsmall}. The methods used there are by no means new: we will essentially reload the classical parametrix machinery for the heat kernel,  carefully track down the lowest order term and the error term, and  estimate the latter uniformly, under suitable global bounds on the curvature tensor. 

 Establishing the sharpness of the exponential constant in \eqref{MS} is not trivial. Ideally, one would like  to fix a point $x$, and then  construct an extremal family of functions in $W^{\a,\na}(M)$ which is equal to $\log (r/d(x,y))$ for large $r$, and $1<d(x,y)<r$.
This is similar in spirit to the standard construction of the so-called ``Moser sequence'' for the classical inequalities (which is  defined on small annuli around a point $x$). However, this  idea  only works for $\alpha=1$, due to the non-smoothness of the distance function on large balls. To overcome this difficulty, we will replace $d(x,y)$ with the function 
 \be b(y)=\bigg(\frac{G_2(x,y)}{\sm^{-1}c_2}\bigg)^{\frac1{2-n}}\sim  d(x,y),\qquad {\text{ as }} d(x,y)\to\infty,\ee
which was featured in several papers by Colding and Colding-Minicozzi (see for ex. \cite{c2}, \cite{cm1}, \cite{cm2}). One should think of $b$ as a smooth, proper replacement of the distance function, on manifolds with EVG. 
The extremal family for \eqref{MS} will then be constructed in Section \ref{sharpness} either by smoothing $\log(r/b)$ on large annuli (which seems to be  a workable method only for  $\alpha=2$), or by taking potentials of suitable compactly supported functions which are essentially equal to $b^{-\a}$ on large annuli, when $2\le \alpha<n/2$ (cf. Remark \ref{b}).

In either case it will be necessary to have at hand sharp estimates for the Green function $G_\a(x,y)$, for large distances from a given point $x$.  In Sections \ref{heatlarge} and \ref{greenlarge}, we will first refine some heat kernel and Green  function estimates for large distances from a fixed point, on manifolds with EVG, due to Li-Tam-Wang \cite{ltw} and Colding-Minicozzi \cite{cm1}. As a result, we will obtain the sharp asymptotic bounds for $\alpha$ even
\be |G_\a(x,y)-\sm^{-1}c_\a d(x,y)^{\a-n}|\le C \wt\Lambda_x\big(d(x,y)\big) d(x,y)^{\a-n},\qquad d(x,y)\ge 1\label{ee1}\ee
\be\bigg(\fint\limits_{r<d(x,y)\le (1+\eta)r}\hskip-1.4em \Big|\nabla_y G_\a(x,y)-\sm^{-1}c_\a \nabla_y d(x,y)^{\a-n}\Big|^2 d\mu(y)\bigg)^{\frac12}\le C_\eta \wt\Lambda_x^{\frac12}(r) r^{\a-n-1},\qquad \eta>0, r\ge1\label{ee2}\ee
where $\wt\Lambda_x(r)$ is a decreasing function such that $\wt\Lambda_x(r)\to0$ as $r\to\infty.$ Such function is obtained from  the volume ratio remainder
\be\Lambda_x(r)=\frac{V_x(r)}{B_nr^n}-\sm\label{vrr}\ee
via  a sort of transform defined in \eqref{Lphi}. In \cite{cm1} estimates \eqref{ee1}, \eqref{ee2} were obtained for $\alpha=2$, but with an unspecified function $\e(r)$ going to 0 as $r\to\infty$, in place of $\wt \Lambda_x(r)$. Using an iteration procedure it would be possible to obtain a similar result for all $\alpha$ even, which  would also be  enough to establish sharpness of \eqref{MS}. We believe, however, that our estimates are of independent interest, since they quantify the  rate of convergence of the heat kernel/Green function in terms of the rate  of convergence of the volume ratio (see Lemma \ref{phitilde} and comments thereafter).

It is noteworthy that in order  to establish the sharpness of the power $\nna$ in the denominator of \eqref{MS}, the bounding  radii of the annuli where the extremal functions are built, need to be finely calibrated in terms of the volume ratio remainder \eqref{vrr}, via its transform $\wt\Lambda_x(r)$ (see \eqref{rrho} and Remark \ref{rrho1}).

 \bigskip

 \section{Background notation}\label{background}

Let $(M,g)$ be a complete, noncompact, $n-$dimensional Riemannian manifold, with metric tensor $g$. The geodesic distance between two points $x,y\in M$ is denoted by $d(x,y)$, and the open geodesic ball centered at $x$ and with radius $r$ will be denoted as $B(x,r)$. The manifold $(M,g)$ is equipped with the natural Riemannian measure $\mu$ which in any chart satisfies $d\mu=\sqrt{|g|} dm$, where $|g|=\det(g_{ij})$, $g=(g_{ij})$ and $dm$ is the Lebesgue measure. 

\smallskip
For each given $p\in M$ denote the tangent space at $p$ by $T_pM$, 
and the unit sphere in such space as $S_p=\{\xi\in T_pM: |\xi|=1\}$ where
\be \lan X,Y\ran=g(X,Y),\qquad |X|=\lan X,X\ran^{1/2},\qquad X,Y\in T_p M.\label{norm}\ee
 \def\ric{{\rm{Ric } }}Covariant derivatives of order $k$ of a smooth function $u$ are denoted as $\nabla^k u$, and their components in a local chart are denoted as $(\nabla^k u)_{j_1...j_k}$. In particular $\nabla^0 u=u$, $(\nabla^1 u)_j=\partial_j u$. Define
 $$|\nabla^k u|^2=g^{i_1j_1}....g^{i_kj_k}(\nabla^k u)_{i_1...i_k}(\nabla^k u)_{j_1...j_k}$$
The gradient and the Laplacian of  a  smooth function $u$ on $M$ are given as  
\be \nabla u=g^{ij}\partial_j u \,\partial_i,\qquad \Delta u=g^{ij}(\nabla^2 u)_{ij}=\frac{1}{\sqrt{|g|}}\,{\partial_i}\left(\sqrt{|g|}\,g^{ij}\partial_j u  \right).\nt\ee
and 
$$|\nabla u|^2=g(\nabla u, \nabla u)=g_{ij}(\nabla u)^i(\nabla u)^j=g^{ij}\partial_i u\, \partial_j u=|\nabla^1 u|^2.$$

Given a smooth  tensor field  $F(x,y)$ on $M\times M$, endowed with the standard product metric, we will use the notation $\nabla_x^1 F$ to denote the covariant derivative of $F$  with respect to $x$. Formally, if $Z$ is any vector field in $T(M\times M)$ identified with $TM\oplus TM$ via the projections $\pi_1(x,y)=x, \pi_2(x,y)=y$, then $(\nabla_x^1)_Z F=\nabla^1_X F$, where $X=\pi_{1*}Z$. Similarly one defines $\nabla_y^1 F.$ Given a smooth function $f(x,y)$,  the second mixed covariant derivative of $f$  will be denoted as $\nabla_x^1\nabla_y^1 f=\nabla_x^1(\nabla_y^1 f)$, i.e. $\nabla_x^1\nabla_y^1 f={\text {Hess} f}\circ (\pi_{1*},\pi_{2*}).$
In the standard coordinate charts of $M\times M$, inherited from $M$, one has that $\nabla_x^1 \nabla_y^1 f=\partial_{x_i,y_j}^2 f dx^i dy^j$, and $|\nabla_x^1\nabla_y^1 f|^2=g^{ij}\partial_{x_i,y_j}^2 f.$

 The curvature tensor in the metric $g$ will be denoted as $R$, the sectional curvature as $K$, and the Ricci curvature as $\ric$.

\begin{rk}\label{bg} Assuming $k^{th}$ order bounded geometry,  it is possible to guarantee that if $0<r_0<\inj(M)$, then  there exists $C_0$, depending on $r_0$ and the curvature bounds,  such that for any $x\in M$, and in any normal  coordinate neighborhood of $x$ we have
\be \big|\partial^j_h g_{ij}\big(\!\exp_x(r\xi)\big)\big|\le C_0,\qquad 0\le j\le k,\; r\le r_0,\ee
where $h$ denotes any multi-index with $|h|=j$ (see for example \cite{e}).

\ni We now discuss briefly the Sobolev spaces. For $u\in C^\infty$ define $\|u\|_{W^{k,q}}=\sum_0^k \Big(\int_M |\nabla^j u|^q d\mu\Big)^{1/q}$ and let $C^{k,q}(M)=\big\{u\in C^\infty(M): \|u\|_{W^{k,q}}<\infty\big\}.$  The Sobolev space $W^{k,q}(M)$ is defined to be the completion of $C^{k,q}(M)$, relative to the norm $\|\cdot \|_{W^{k,q}}$, and it can be  viewed as a subspace of $L^q(M)$. Endowed with the obvious norm, still denoted as $\|\cdot \|_{W^{k,q}}$, the space $W^{k,q}$ is a Banach space. The Sobolev space $W_0^{k,q}(M)$ is the closure of $C_c^\infty(M)$ in $W^{k,q}(M)$.
We have already defined the Bessel potential spaces $H^{\a,m}(M)$ in the introduction. We remark here that on manifolds with $\ric\ge0$ the Riesz transform $\nabla(-\Delta)^{\frac12}$ is bounded on $L^p$ for $1<p<\infty$, and $\|\nabla u\|_p$ is equivalent to  $\|(-\Delta)^{\frac12}u\|_p$, due to results in \cite{acdh}, Thm 1.4, since such $M$ satisfies the doubling property and the heat kernel estimates $H(t,x,y)\le C/V_x(\sqrt t),\; |\nabla H(t,x,y)|\le C t^{-1/2}/V_x(\sqrt t)$  (see Section 3). Hence we can use $\|u\|_p+\|D^\a u\|_p$ as a norm in $H^{\a,p}$ and we also have $W_0^{\a,p}\subseteq W^{\a,p}\subseteq H^{\a,p}$, for $p>1$, $\a\in \N$. We will not be concerned  with the issues related to when equality holds in the previous inclusions. Suffices to say that it's certainly true that for $p\ge1$ we have $W_0^{\a,p}=W^{\a,p}$ when $\alpha=1$, and for $\alpha\ge 2$ provided that $\inj(M)>0$ and  $|\nabla^j \ric|\le C$, for $j\le \a-2$ (see \cite{h}, Thms. 2.7, 2.8, and \cite{ve} for related counterexamples). Regarding the spaces $H^{\a,p}$, it's clear that for $p>1$  we have $W_0^{\a,p}=W^{\a,p}=H^{\a,p}$ when $\a=1$, and for higher values of $\alpha$ it is believed to be true again when $\inj(M)>0$ and $|\nabla^j \ric|\le C$, $j\le \a-2$ (see [GP] for results in this direction when $\alpha=2$), and it's certainly true for all $\alpha$ if $\inj(M)>0$ and $|\nabla^j R|\le C$ for all $j$ (\cite{tr}, p.320). 

The cut locus of $x$ can be described as $\cut_x=\exp_x\{c(\xi)\xi,\;\xi\in S_x\}$, where $\exp_x$ is the exponential map at $x$, and $c(\xi)<\infty$ is the distance from $x$ to its cut point. Then, 
\be \exp_x:\;D_x:=\{t\xi: 0< t<c(\xi),\; \xi\in S_x\} \to M_x:=M\setminus (\cut_x\cup\{x\})\label{DpMp}\ee
is a diffeomorphism. Moreover, $\cut_x$ has zero measure and for any $f$ integrable on $M$ we can write the integral of $f$ in geodesic polar coordinates as 
 \be\int_M f d\mu=\int_{M_x} fd\mu = \int_{S_x}d\mu_x(\xi) \int_0^{c(\xi)} f(\exp_x t\xi)\sqrt{|g|}(t\xi) t^{n-1}dt,\label{int}\ee
where $d\mu_x(\xi)$ is the measure on $S_x$ induced by the Lebesgue measure (see \cite{cha}, \textsection III.3]).

In this notation we have, for each $x\in  M$, 
\be V_x(r)=\int_{B(x,r)} d\mu =\int_0^r A_x(s) ds\ee
where 
\be A_x(r)=\int_{\{r<c(\xi)\}}\sqrt{|g|}(r\xi)r^{n-1} d\mu_x(\xi)\ee
which coincides with the $(n-1)$-dimensional Hausdorff measure of $\partial B(x,r)$, for a.e. $r>0$.

By the Bishop comparison theorem the functions $A_x(r)/r^{n-1}$ and $V_x(r)/r^{n}$ are both decreasing for $r>0$, and
\be\lim_{r\to0} \frac{V_x(r)}{B_n r^n}=1.\ee
Using that $B(x,r)\subseteq B(y,r+d(x,y))$ one sees that
\be \sm:=\lim_{r\to+\infty}  \frac{V_x(r)}{B_nr^n} \label {sigma}\ee
is independent of $x$, from which we deduce
\be \sm B_nr^n\le V_x(r)\le B_n r^n,\qquad x\in M,\;\; r>0.\label{Evg3}\ee
Similarly, since $A_x(r)/r^{n-1}$ is decreasing, \eqref{sigma} and \eqref{Evg3} imply (and are implied by)
\be \sm=\lim_{r\to+\infty}  \frac{A_x(r)}{nB_n r^{n-1}},\qquad  \sm nB_n r^{n-1}\le A_x(r)\le nB_n r^{n-1}.\label {sigma4}\ee

Clearly, $\sm\le1$ and it turns out, from the Bishop comparison theorem,  that if $\sm=1$ then $(M,g)$ is isometric to $\R^n$. As a consequence of this, and of the  proof of the Bishop comparison, the function $V_x(r)/r^n$ is strictly decreasing, unless $(M,g)$ is isometric to $\R^n$.

\bigskip

\section{Uniform heat kernel estimates}

In this section  $(M,g)$ will denote  a complete, $n-$dimensional Riemannian manifold with nonnegative Ricci curvature.  The unique positive, symmetric, minimal heat kernel, i.e. the fundamental solution of $\Delta-\partial_t=0$ on $M$, will be denoted as $H(x,y,t).$ For the existence of $H$ and its basic properties see for example \cite{li} or \cite{st}.

\subsection{Uniform global bounds.} In this section we recall some known Gaussian bounds for the heat kernel and its space and time derivatives. First, we have the well-known bounds due to Li-Yau \cite{ly}: for any $\delta\in (0,1)$, there exist $\Cl{1}, \Cl{2}>0$, depending only on $\delta$ and $n$, such that for any $x,y\in M$ and any $t>0$
\be \Cr{1} V_x(\sqrt t)^{-1} \exp\left(-\frac{d^2(x,y)}{4(1-\delta)t}\right)\le H(x,y,t)\le \Cr{2} V_x(\sqrt t)^{-1} \exp\left(-\frac{d^2(x,y)}{4(1+\delta)t}\right).\label{hrough}\ee

Regarding the global gradient estimate,  there exists $\Cl{3}>0$ depending only on $n$  such that for all $x, y\in M$, and all $t>0$
\be |\nabla_y H(x,y,t)|\le \Cr{3}t^{-\frac1 2}\left(1+\frac{d^2(x,y)}t\right) H(x,y,t).\label{gradh}\ee
 Estimate \eqref{gradh} is due to Souplet-Zhang \cite{sz}.
 
 \medskip
 The following time derivative estimate was obtained by Grigor'yan \cite{g}, Cor. 3.1]:   there exists $\Cl{gr}$ such that for all $x, y\in M$, and all $t>0$
 \be  \left|\frac{\partial H}{\partial t}(x,y,t)\right|\le \Cr{gr}  t^{-1} V_x(\sqrt t)^{-1}\bigg(1+\frac{d^2(x,y)}{4t}\bigg)^{2+3n/4} \exp\bigg(-\frac{d^2(x,y)}{4t}\bigg).\label{grig}\ee
 
 Finally, a uniform bound on the mixed second  covariant derivative was obtained by Li Jiayu \cite{lij}: if  $\alpha>1$, then there exists $\Cl{li}$ depending on $n$ and $\alpha$ such that for all $x,y\in M$ and all $t>0$
  \be |\nabla_x^1\nabla_y^1 H(x,y,t)|\le \Cr{li} t^{-1} H(x,y,t)+3\alpha \frac{\p H}{\p t}(x,y,t)-\alpha\frac{|\nabla_y H (x,y,t)|^2}{H(x,y,t)}. \label{li}\ee  
 
 Under the assumptions $\ric\ge0$ the Bishop comparison gives $V_x(r)\le  |B_{\sR^n}(r)|$ for all $r>0$.  Assuming    $\inj(M)>0$, and sectional curvature $K\le k_0$, some $k_0>0$, (in particular under $1^{st}$ order bounded geometry),  the G\"unther volume comparison theorem   (\cite{lee},  Thm. 11.14)
 gives 
 \be  V_x(r)\ge  \vol_{k_0}(r),\qquad \forall x\in M,\; 0<r<\min\Big\{\inj(M),\frac\pi{\sqrt{k_0}}\Big\}.\label{Gun0}\ee
where $\vol_{k_0}(r)=\omega_{n-1}\ds\int_0^r\Big[\frac1{\sqrt{k_0}}\sin\big(\sqrt {k_0}\,t\big)\Big]^{n-1}dt$ denotes the volume of the geodesic ball of radius $r$ in the space form of constant sectional curvature $k_0$ (sphere).
Therefore, if $r_1=\min\big\{\inj(M),\pi/\sqrt{k_0}\big\}$, then the Bishop comparison implies  that for any $\rb>0$, $0<r\le \rb$ we have  \be \frac{V_x(r)}{ |B_{\sR^n}(r)|}\ge\frac {V_x(\rb)}{ |B_{\sR^n}(\rb)|}\ge \frac 1{ |B_{\sR^n}(\rb)|} \min\big\{\vol_{k_0}(\rb),\vol_{k_0}(r_1)\big\}\ee and
 \be  V_x(r)\ge \Cl{gun}r^n,\qquad 0<r\le \rb\label{Gun}\ee
some  $\Cr{gun}>0$ independent of $x$ (and depending on $\rb$, in general).

 \smallskip
 Let
 \be E(x,y,t)=(4\pi t)^{-{\frac n 2}}\exp\left(-\frac{d^2(x,y)}{4t}\right).\label{E}\ee
Taking $\delta=\frac1 2$ in  the bounds  \eqref{hrough}, \eqref{gradh} and  using \eqref{Gun}  together with the Bishop comparison,   we obtain that for all $x,y\in M$ and for  $0<t\le 1$ 
 \be \Cl{7} E(x,y,t/2)\le H(x,y,t)\le \Cl{8} E(x,y,3t/2)\label{hrough1}\ee
 \be |\nabla_y H(x,y,t)|\le\Cl{9} t^{-\frac1 2} E(x,y,3t) \label{gradh1}\ee
 for some constants $\Cr{7},\Cr{8},\Cr{9}$ independent of $x,y,t$, where the left inequality in \eqref{hrough1} actually holds for all $t>0$. Clearly, $|\nabla_x E|$ satisfies the same bound as in \eqref{gradh1}, but within the domain of differentiability of $d(x,y)$ (defined in \eqref{DpMp}), in particular a.e.:
  \be |\nabla_y E(x,y,t)|\le\Cl{9} t^{-\frac1 2} E(x,y,3t),\qquad y\in M_p.  \label{gradE1}\ee
 
 \begin{rk}\label{rem1} Note that \eqref{gradh1},    and the right inequality in \eqref{hrough1}, hold  for $0<t\le T$, any $T>0$, with $\Cr8$ and $\Cr9$ depending on $T$, and this is true under the assumption $\ric\ge0$, $\inj(M)>0$,  and $K$ bounded.    Assuming Euclidean volume growth instead of $K$ bounded, such inequalities hold for all $t>0$, since $V_x(r)\ge \sm |B_{\sR^n}(1)|\,r^n.$
 \end{rk}
 
 Later we will also need the following Lipschitz estimates for $H$: assuming $K$ bounded, $T, \eta>0$, there are $\Cl{Hol1}, \Cl{Hol2}$ such that for any $R>(1+\eta)r$, $x'\in B(x,r)$,  $y\notin B(x,R)$, and for all $0<t\le T$
 \be |H(x,y,t)-H(x',y,t)|\le \Cr{Hol1} d(x,x') t^{-\frac 12}V_x(\sqrt t)^{-1}\exp\bigg(\!\!-\frac{d^2(x,y)}{c_\eta t}\bigg)\label{hol1}\ee
\be  |\nabla_y H(x,y,t)-\nabla_y H(x',y,t)|\le \Cr{Hol2} d(x,x') t^{-1}V_x(\sqrt t)^{-1}\exp\bigg(\!\!-\frac{d^2(x,y)}{c_\eta t}\bigg)
 \label{hol2}\ee
 for some $c_\eta>0$ (specifically, we can take $c_\eta=6(1+\eta)^2/\eta^2$).  Estimate \eqref{hol1} follows easily from the mean value theorem along a minimizing geodesic joining $x$ and $x'$, and  the gradient estimate \eqref{gradh}, in the $x$ variable, followed by  \eqref{hrough} (with $\delta=1/2$).  Estimate \eqref{hol2} follows in a similar manner, but using instead  the mixed gradient estimate \eqref{li}, followed by   \eqref{grig} (with $\a=2$), and \eqref{hrough}.

 \subsection{Uniform asymptotic bounds: small times,  small distances, bounded geometry.} \label{heatsmall}
 
 \medskip\medskip
 It is well known that locally the heat kernel on $M$, for small times, behaves like the Euclidean heat kernel, i.e. in a suitable sense
 \be H(x,y,t)\approx (4\pi t)^{-{\frac n 2}}\exp\left(-\frac{d^2(x,y)}{4t}\right),\ee
 for small $d(x,y)$ and small $t$ (cf. Remark \ref{asymp}).
 In the following proposition  we quantify this statement under the assumption of  $1^{st}$  order bounded geometry, and give a corresponding statement for the gradient of $H$, under  $2^{nd}$  order bounded geometry. While the results seem somewhat intuitive, we were not able to find them in the vast  literature, in the particular form stated below. Once again, the main issue is to guarantee uniformity in the error term. 
 
 Once and for all let us fix an $r_0$ such that 
 \be 0<r_0<\inj(M).\ee
  \begin{prop}\label{smalld} If $M$ has $1^{st}$ order bounded geometry,  then there exists $\Cl{4}>0$ depending only on  $n$, $r_0$,  and the curvature bounds such that for $0<t\le 1$ and $x,y\in M$
\be \left| H(x,y,t)-   (4\pi t)^{-{\frac n 2}}\exp\left(-\frac{d^2(x,y)}{4t}\right)\right|\le \Cr{4} t^{-{\frac n 2}+\frac1 2} \exp\left(-\frac{d^2(x,y)}{24t}\right).
\label{hsmall}\ee 
If $M$ has $2^{nd}$ order bounded geometry then there exists $\Cl{5}>0$ depending only on  $n$, $r_0$,  and the curvature bounds such that for $0<t\le 1$, $x\in M$, and $y\in M_x$
\be \left| \nabla_yH(x,y,t)-  \nabla_y (4\pi t)^{-{\frac n 2}}\exp\left(-\frac{d^2(x,y)}{4t}\right)\right|\le \Cr{5} t^{-{\frac n 2}} \exp\left(-\frac{d^2(x,y)}{24t}\right).
\label{gradhsmall}\ee

 \end{prop}
 
 \begin{rk}\label{asymp} Note that \eqref{hsmall} implies that within each region of type $d^2(x,y)\le Ct$, some fixed $C>0$, we have
 \be H(x,y,t)\sim (4\pi t)^{-{\frac n 2}}\exp\left(-\frac{d^2(x,y)}{4t}\right),\qquad t\to0,\ee
 where $f(t)\sim g(t)$ means that $f(t)/g(t)\to1$, as $t\to0$.
 
 \smallskip
 \ni{\bf Proof.}   
  We note first that \eqref{hsmall}, \eqref{gradhsmall}  are only meaningful for $d(x,y)$ small. Indeed, assuming $d(x,y)\ge a>0$ (some fixed $a>0$), one can easily prove both bounds using the triangle inequality and  \eqref{hrough1}, \eqref{gradh1}.  
Through out the proof we will then assume $d(x,y)\le r_0$. 

Let us start with proving  \eqref{hsmall}. We combine the usual construction of the parametrix for  heat kernel with the global Gaussian bounds in \eqref{hrough}, \eqref{gradh}.
In geodesic normal coordinates at $x$ let
 \be u_0(x,y)=|g|^{-\frac 1 4}(r\xi)\ee and let $\phi\in C^\infty(M)$ be so that $\phi=1$ on $B(x,r_0/2)$, $\phi=0$ outside $B(x,r_0)$, $0\le \phi\le 1$ and $|\nabla\phi|\le 2/r_0.$
 
 Letting
 \be F(x,y,t)=\phi(y)u_0(x,y)E(x,y,t)\ee
 by standard calculations (see \cite{li}, pp. 109-111) we have
 \be H(x,y,t)-F(x,y,t)=\int_0^t\int_M H(z,y,t-s)\big(\Delta_z -\partial_s\big) F(x,z,s)dzds,\label{FH1}\ee
 and for $z\in B(x,r_0)$ 
 \be \big(\Delta_z -\partial_t\big) F(x,z,t)=\phi(\Delta_z u_0) E +u_0 E \Delta_z\phi+2
\lan \nabla_z\phi,\nabla_z (u_0 E)\ran.\label{FH2}\ee

Inserting \eqref{FH2} into \eqref{FH1}, using Green's identity in the integrals corresponding to the first two terms  in \eqref{FH2}, we obtain
\be\ba H(x,y,t)-F(x,y,t)&=-\int_0^t\int_M \lan \nabla_z H(z,y,t-s),\nabla_z(\phi u_0)\ran E(x,z,s)dxds\cr&
+\int_0^t\int_M H(z,y,t-s)\lan\nabla_z E(x,z,s),u_0\nabla_z\phi-\phi\nabla_z u_0\ran dxds.\ea
 \label{formula}\ee
Using the bounds \eqref{hrough1}, \eqref{gradh1} we obtain, for $0<t\le1$, $y\in B(x,r_0),$
\be\ba |&H(x,y,t)-F(x,y,t)|\le\cr&
\le \Cl{10}\!\int_0^t\big[(t-s)^{-\frac1 2}+s^{-\frac1 2}\big]ds\int_{B(x,r_0)} \!\!(u_0+|\nabla_zu_0|)(x,z)E(z,y,3(t-s))E(x,z,3s)dz\label{HF1}\ea\ee
and using the hypothesis of bounded geometry \eqref{M4} (see Remark \ref{bg}), the lower bound for $H$ in \eqref{hrough1} (valid for all $t>0$, see Remark \ref{rem1}), and the semigroup property 
\be\ba  |H(x,y,t)&-F(x,y,t)|\le\Cl{11}\!\int_0^t s^{-\frac1 2}ds\int_{M} H(z,y,6(t-s))H(x,z,6s)dz\cr&
= \Cr{11} H(x,y,6t)\int_0^t s^{-\frac1 2} ds\le \Cl{12} t^{\frac1 2-{\frac n 2}} \exp\left(-\frac{d^2(x,y)}{24t}\right).\label {HF2}\ea\ee

Finally, note that $u_0(x,x)=\phi(x)=1$, hence using \eqref{M4} we get $|u_0(x,y)\phi(y)-1|\le\Cl{13} d(x,y)$, for $y\in B(x,r_0)$, and therefore
\be|F(x,y,t)-E(x,y,t)|= |E(x,y,t)||u_0(x,y)\phi(y)-1|\le \Cl{14}t^{\frac1 2-{\frac n 2}}  \exp\left(-\frac{d^2(x,y)}{8t}\right).\label{FE}\ee

To prove \eqref{gradhsmall},  we use again  the representation formula \eqref{FH1}. 
This time after inserting \eqref{FH2} into \eqref{FH1} we take directly the gradient in the $y$ variable, obtaining
\be \ba|\nabla_y H(x,y,t)&-\nabla_y F(x,y,t)|\le \cr& \le \C\int_0^t\int_M|\nabla_y H(z,y,t-s)|\,|E(x,z,s)|\big(u_0+|\nabla_z u_0|+|\Delta_z u_0|\big)(x,z)dzds\cr& \hskip3em +\C\int_0^t\int_M|\nabla_y H(z,y,t-s)|\left|\nabla_z E(x,z,s)\right|u_0(x,z)dzds.
\ea\ee
Assuming  $2^{nd}$ order bounded geometry, using the bounds \eqref{hrough1}, \eqref{gradh1}, and arguing just as in \eqref{HF1}, \eqref{HF2}  we obtain, for $0<t\le1$, $y\in B(x,r_0),$
\be\ba|\nabla_y H(x,y,t)-\nabla_y F(x,y,t)|&\le\C \int_0^t (t-s)^{-\frac1 2}s^{-\frac1 2}ds\int_{B(x,r_0)}E(z,y,3(t-s))E(x,z,3s)dz\cr&\le \C t^{-{\frac n 2}} \exp\left(-\frac{d^2(x,y)}{24t}\right).\ea\ee
  Finally,
  \be\ba|\nabla_y F-\nabla _y E|&=\left|\nabla_y\big(E(\phi u_0-1)\big)\right|\le|\nabla_y u_0| E+|\nabla_y E| |\phi u_0-1|\cr& \le \C t^{-{\frac n 2}} \exp\left(-\frac{d^2(x,y)}{8t}\right)\ea\ee
  arguing as in \eqref{FE}.

\qed

\subsection{Uniform asymptotic bounds: large distances from a fixed point}\label{heatlarge}

\bigskip
Let $M$ have Euclidean volume growth i.e.
\be \sx(r):=\frac{V_x(r)}{B_n r^n}\searrow\sm>0.\label{v1}\ee 
In \cite{ltw}, Li-Tam-Wang proved the following heat kernel bound: there is $\Cl{15}>0$ depending only on $n$ and $\sm$, such that for all $x,y\in M$, $t>0$ and $\delta>0$
\be\ba \frac{(4\pi t)^{-{\frac n 2}}}{\tx\big(\delta d(x,y)\big)}&\exp\left(-\frac{1+9\delta}{4t} d^2(x,y)\right)\le H(x,y,t)\le\cr&\le \big(1+\Cr{15}(\delta+\beta(x,y,\delta)\big)\frac{(4\pi t)^{-{\frac n 2}}}{\sm}\exp\left(-\frac{1-\delta}{4t} d^2(x,y)\right)\ea\label{LTW}\ee
where
\be\beta(x,y,\delta)=\delta^{-2n}\max_{r\ge (1-\delta)d(x,y)}\left(1-\frac{\tx(r)}{\tx(\delta^{2n+1}r)}\right),\ee
and where
\be\tx(r)=\frac{A_x(r)}{n B_n r^{n-1}}.\ee
(Note that in  the notation used in \cite{ltw},  $\tx(r)=A_x(r)/(nr^{n-1}).$)

In particular, since for fixed $\delta>0$ we have $\beta(x,y,\delta)\to0$, as $d(x,y)\to+
\infty$, then \eqref{LTW} gives, within each region where $d^2(x,y)\le C t$ (some fixed $C>0$)
\be H(x,y,t)\sim\frac{1}{\sm} (4\pi t)^{-{\frac n 2}}\exp\left(-\frac{d^2(x,y)}{4t}\right),\qquad d(x,y)\to +\infty\ee
(cf. Remark \ref{asymp}).
In the following proposition we work out \eqref{LTW} a bit more explicitly, in terms of the \emph{volume ratio remainder} at a given point $x\in M$, defined as
\be \L_x(r)=\sx(r)-\sm=\frac{V_x(r)}{B_nr^n}-\sm.\label{R}\ee

Let us first note that using $\tx(r)$ in place of $\sx(r)$ in \eqref{R} gives an equivalent remainder:

\begin{lemma}\label{txsx} There is $C_M\in (0,1]$ such that 
\be C_M\big(\sx(r)-\sm\big)\le \tx(r)-\sm\le \sx(r)-\sm.\ee
\end{lemma} 
\smallskip
\ni{\bf Proof.} On one hand we have
\be \sx(r)=\frac1{B_n r^n}\int_0^r A_x(t)dt=\frac1{n r^n}\int_0^r \tx(t) n t^{n-1}dt\ge  \frac{\tx(r)}{n r^n}\int_0^r  n t^{n-1}dt=\tx(r).\ee

On the other hand, assuming $\sm<1$ (if $\sm=1$ then $\sx(r)=\tx(r)=1$ for all $r$), the isoperimetric inequality \eqref{isop} gives
\be \tx(r)-\sm\ge \sm^{\frac 1n} \sx(r)^{\frac{n-1}n}-\sm=\sm\bigg(\Big(\frac{\sx(r)}\sm\Big)^{\frac{n-1}n}-1\bigg)
\ge \psi(\sm^{-1})\big(\sx(r)-\sm\big)
\ee
where $\psi(t)=(t^{\frac{n-1}n}-1)/(t-1)$, a decreasing function on $(1,\infty)$.

\qed

\smallskip

\begin{de}\label{tilde}
Given $\phi:[0,\infty)\to(0,\infty)$, differentiable,  decreasing, and $\phi(r)\to0$ as $r\to\infty$, let, for $\delta,r>0$, 
\be \wt\phi(r)=\min_{\delta>0} \big(\delta+\delta^{-2n}\phi(\delta^{2n+1}r)\big).\label{Lphi}\ee\end{de}
It's clear that $\wt\phi(r)$ is well defined for all $r>0$, and decreasing. In  Lemma \ref{phitilde} we will  show that  $\wt\phi(r)\to0$ as $r\to+\infty
$, and in fact we will give an effective, and  essentially optimal, bound for $\wt\phi$ in terms of $\phi$. 

 \begin{prop}\label{largedH} If $M$ has EVG, then there exists $\Cl{16}$ depending only on  $n,B,\sm$  such that for any $x,y\in M$  with  $d(x,y)\ge1$, and for any $t>0$
\be \Big| H(x,y,t)- \frac{1}{\sm}  E(x,y,t)\Big|\le \Cr{16} \;\wt\L_x\big(d(x,y)\big)\, E(x,y,3t).
\label{hlarge}\ee

 \end{prop}
 
 \ni{\bf Proof.} For $0<\delta\le \half$ we have
 \be\ba\beta(x,y,\delta)&\le\delta^{-2n}\max_{r\ge \frac12 d(x,y)}\left(\frac{\tx(\delta^{2n+1}r)-\sm}{\tx(\delta^{2n+1}r)}-\frac{\tx(r)-\sm}{\tx(\delta^{2n+1}r)}\right) \\&\le \delta^{-2n}\max_{ r\ge \frac12 d(x,y)}\frac{\tx(\delta^{2n+1}r)-\sm}{\tx(\delta^{2n+1}r)}\le \delta^{-2n}\max_{ r\ge \frac12 d(x,y)}\frac{\tx(\delta^{2n+1}r)-\sm}{\sm}\\&=\frac{\delta^{-2n}}{\sm}\Big(\tx\big(\half\delta^{2n+1}d(x,y)\big)-\sm\Big).\ea\label{beta}\ee
Also, 
\be\frac1\sm-\frac1{\tx\big(\delta d(x,y)\big)}\le\frac{\tx\big(\delta d(x,y)\big)-\sm}{\sm^2}\le2 \frac{\delta^{-2n}}{\sm^2}\Big(\tx\big(\half\delta^{2n+1}d(x,y)\big)-\sm\Big)\ee
since $\tx(\lambda r)/\lambda$ is decreasing in $\lambda.$

Since 
 \be|e^{-\lambda_1x}-e^{-\lambda_2 x}|\le 3|\lambda_1-\lambda_2|e^{-x/3},\qquad \lambda_1,\lambda_2\ge \half,\; x>0.\label{exp}\ee
and changing $\delta$ into $2^{\frac1 {2n+1}}\delta$ and using \eqref{LTW}  we  get, for $\delta\le\quarter$, 
 \be \Big|H(x,y,t)-\frac1\sm E(x,y,t)\Big|\le \C\Big[\delta+\delta^{-2n}\Big(\tx\big(\delta^{2n+1}d(x,y)\big)-\sm\Big)\Big]\; E(x,y,3t).\ee
\ni It's easy to check that the same estimate holds for $\delta>\quarter$, so, using Lemma \ref{txsx},  we get \eqref{hlarge}.

\qed

The next Lemma provides more explicit information about the function $\wt\phi$:

\begin{lemma} \label{phitilde} Let $\phi,\wt\phi$ be as in Definition \ref{tilde}, and let $t(r)$ be the unique solution of $t/\phi(t)=r$. Then $t(r)$is differentiable, strictly increasing to $+\infty$, and
\be \phi(t(r))^{\frac1{2n+1}}\le \wt\phi(r)\le 2\phi(t(r))^{\frac1{2n+1}},\qquad r>0.\label{tilde1}\ee
\end{lemma}

\ni{\bf Proof.}   The statement about $t(r)$ is clear since $t/\phi(t)$ is strictly increasing, and $\phi(t)\to0$, as $t\to\infty.$  Let  $\delta_r^{2n+1}=t(r)/r=\phi(t(r))$, which is equivalent to $\delta_r=\delta_r^{-2n}\phi(\delta_r^{2n+1}r)$. For fixed $r$ let $F(\delta)=\delta+\delta^{-2n}\phi(\delta^{2n+1}r)$, so that $F(\delta_r)=2\delta_r$. If the minimum of $F$ is attained at $\delta_{r,0}$, then an easy calculation shows that 
$ F(\delta_{r,0})-\delta_{r,0}= \delta_{r,0}^{-2n}\phi(\delta_{r,0}^{2n+1}r)\le\frac1{2n}\delta_{r,0}$, and since $F(\delta)-\delta$ is decreasing, this implies that $\delta_r\le \delta_{r,0}$. Hence, we obtain
\be \delta_r\le \delta_{r,0}\le F(\delta_{r,0})\le F(\delta_r)=2\delta_r\ee
which implies \eqref{tilde1}.

\qed

In practical applications, given an estimate $\Lambda_x(r)\le \phi(r)$, for $r\ge r_0$ with $\phi$ differentiable, decreasing, and going to 0 at infinity, the above lemma gives the  rather explicit  upper bound $\wt\Lambda_x(r)\le 2\phi(t(r))^{\frac1{2n+1}}$, for $r\ge r_0$.

For example, when $\phi(r)= C r^{-a}(\log r)^{-b}$ with $a\ge0$, $b\in\R$, and $b>0$ if $a=0$, then one can check (using the Lambert function when $b\neq 0$) that for large $r$
\be t(r)\asymp  r^{\frac 1{1+a}}(\log r)^{-\frac b{1+a}} \qquad  \wt\phi(r)\asymp   r^{-\frac a{(1+a)(2n+1)}}(\log r)^{-\frac b{(1+a)(2n+1)}}.\label{Hab}\ee
where  $f(r)\asymp g(r)$ means $c_1 g(r)\le f(r)\le c_2 g(r)$ for large $r$.


Asymptotically Locally Euclidean (ALE) manifolds provide an  important family of manifolds with EVG such that  $\Lambda_x(r)\le C r^{-a}$ for some $a>0$, some $x \in M$, and for large $r$. For the definitions and properties of ALE manifolds see for example \cite{bkn}, \cite{afm}.

We point out, however, that it is relatively easy  to construct manifolds with EVG, such that $\Lambda_x(r)$ has arbitrarily slow rate of decay at infinity. For example, one can do so by considering rotationally invariant metrics on $\R^n$.

\bigskip

We will now derive sharp asymptotic estimates for the gradient of the heat kernel, and for large distances. In general we do not expect optimal pointwise bounds (see Remark \ref{rem3}) but  sharp, aymptotic, $L^2$  integral bounds on annuli will suffice for our purpose.

\smallskip
For any $r,R>0$ and $x\in M$ let 
\be A_{r,R}(x)=\{y\in M: r\le d(x,y)\le R\}\ee
for simplicity we will often omit the dependence on $x$ in what follows.

The integral average over a Borel measurable set $A$ will be denoted as 
\be \fint_A =\frac{1}{\mu(A)}\int_A.\ee

We will also use the notation
\be E(r,t)=(4\pi t)^{-{\frac n 2}} \exp\bigg(\!\!-\frac{r^2}{4t}\bigg)\ee
so that $E(x,y,t)=E\big(d(x,y),t\big)$.
\medskip
\begin{prop}\label{largedgradH} If $M$ has EVG
 then there exists $\Cl{17}$ depending only on  $n,\sm$  such that for any $x\in M$,  $\eta>0$, $t>0$, and   $r\ge1$ we have
 \be \bigg(\fint_{A_{r,(1+\eta)r}} \Big|\nabla_y H(x,y,t)- \frac{1}{\sm}\nabla_y  E(x,y,t)\Big|^2 dy\bigg)^{\frac1 2}\!\!\le \Cr{17}\Big(1+\frac{1}{\sqrt\eta}\Big)
\wt\L_x(r)^{\frac12} r^{-\frac12}t^{-\frac1 4}E(r,3t).\label{L2estH}
 \ee
 
 \end{prop}

\begin{rk}\label{nablaE1}  In \eqref{L2estH} the vector field $\nabla_y  E(x,y,t)$ is intended to be defined for  a.e. $y$, in particular on $M_x=M\setminus\cut_x\setminus \{x\}.$ The same type of remark can be made regarding the integral in \eqref{gradgsharp1}.\end{rk}

\smallskip
 \ni{\bf Proof.}  Since we are planning on using Stokes' theorem we need to regularize $E(d(x,y),t)$, as a function of $y$. Using a mollifier it's possible to show that given any compact set $K\subseteq M$ and $x\notin K$ then there exists a sequence of smooth, compactly supported functions $\{d_j \}$ (depending on $x$ and $K$) such that  for any $y\in K$ and any $j$
 
 \smallskip

 i) $\; |d_j(x,y)-d(x,y)|<\dfrac 1 j$
 
 \smallskip
ii) $\; |\nabla_y d_j(x,y)-\nabla_y d(x,y)|\to0$  as   $j\to\infty$, if $K\subseteq M\setminus \cut_x$

\smallskip

iii) $\;|\nabla_y d_j(x,y)|\le 2$ 

\smallskip
 
iv) $\; \Delta_y d_j(x,y)\le \dfrac{n-1} {d(x,y)} +\dfrac 1 j$ 
 
 \medskip
\ni (cf. \cite{ltw}, \cite{c1}). The last estimate is basically an approximated  form of the Laplacian comparison theorem, i.e.  $\Delta_y d(x,y)\le(n-1)/d(x,y)$, which holds in the sense of distributions and also inside $M_x$.

Using these estimates, it's straightforward to check that for any $y\in K$ and $t>0$ and for any $j\ge j_0(x,K)$ large enough 
 
 \be |E(d_j(x,y),t)-E(d(x,y),t)|\le \frac{\Cl{18}}j \; \dist(x,K)^{-1-n}\label{reg1}\ee

  \be\qquad \lim_{j\to\infty} |\nabla_y E(d_j(x,y),t)-\nabla_y E(d(x,y),t)|=0, \quad {\text{ if }} K\subseteq M\setminus \cut_x\label{reg2}\ee
  
   \be |\nabla_y E(d_j(x,y),t)|\le {\Cl{19}} \; \dist(x,K)^{-1-n}\label{reg3}\ee
 
 \be \Delta_y E(d_j(x,y),t) \ge-\frac n 2 t^{-{\frac n 2}-1}-\frac{\Cl{20}}j \big(\dist(x,K)^{-1-n}+\dist(x,K)^{-3-n}\big)\label{reg4}\ee
 
 where the constants $\Cr{18}, \Cr{19}, \Cr{20}$ depend only on $n$.
 \smallskip
 Note that the Laplacian comparison (in the distributional sense) gives the estimate
 \be \Delta_y E(d(x,y),t) \ge\partial_t E(d(x,y),t)\ge-\frac n 2 t^{-{\frac n 2}-1}.\ee

 For simplicity, in the estimates below  we will mute the variables $x,y$ in the integrands, and we will let 
 $$E_j=E_j(x,y,t)=E(d_j(x,y),t).$$

  Since $\cut_x$ has zero measure, then for almost all $r>0$ the set $\partial B(x,r)\cap\cut_x$ has zero $(n-1)-$dimensional Hausdorff measure. Following \cite{ltw} we shall call such $r$ ``permissible''. If  $r,R$ are permissible, with $R>r>1$ , then $\Ar$ is ``almost $C^1$'', in the sense that around each point of $\partial \Ar\setminus \cut_x$ the boundary of $\Ar$ can be represented by a $C^1$ function. For such domains Stokes' theorem works (see for example \cite{m}, Thm. 9.6, for the $\R^n$  version, which can be extended to Riemannian manifolds). Hence, for all permissible $r,R$ with $R>r>1$ we have
  
 \be\ba &\mathop\int\limits_{\Ar\setminus\cut_x} \big|\nabla H-\smi\nabla E|^2d\mu \le \liminf_j\mathop\int\limits_{\Ar\setminus\cut_x} \big|\nabla H-\smi\nabla E_j|^2d\mu\cr&= \liminf_j\int_{\Ar} \big|\nabla H-\smi\nabla E_j|^2d\mu\cr&=\liminf_j\bigg(\!\!\!-\!\!\int_\Ar \!\!(H-\smi E_j)\Delta(H-\smi E_j)d\mu+\int_{\p \Ar} \!\!(H-\smi E_j)\nabla (H-\smi E_j)\cdot \bN d\mut\bigg)\ea\label{fatoustokes}\ee
 
\ni where $\bN$ is  the outward unit normal vector field along $\p \Ar$ and $\mut$ the induced measure on $\p\Ar$, and where the first inequality follows from \eqref{reg2} and Fatou's lemma.

\smallskip
\ni To estimate the first integral note that from \eqref{hlarge} and \eqref{reg1} we have (with $d=d(x,y)\ge r\ge1$)
\be |H-\smi E_j|\le \Cl{21}\, \wt\L_x\big(d(x,y)\big)\, E(d,3t)+\frac{\Cl{22}} j\le \Cr{21}\, \wt\L_x(r)\, E(r,3t)+\frac{\Cr{22}} j.\label{i1} \ee
Using the heat equation and the time derivative estimate \eqref{grig} we have
\be|\Delta H|=|\partial_t H|\le \C\, t^{-{\frac n 2}-1}.\label{i2}\ee

To estimate the integral of $|\Delta E_j|$ over the annulus, we use the same idea as in \cite{ltw}, proof of Thm. 1.1, \cite{c1}, proof of Lemma 1.4. From the Laplacian comparison \eqref{reg4} we have
\be\Delta E_j\ge - \frac n 2\,t^{-{\frac n 2}-1}-\frac {\Cl{23}}j\,\qquad \big(\Delta E_j\big)_-\le  \frac n 2\,t^{-{\frac n 2}-1}+\frac {\Cr{23}}j
\ee

where $\big(\Delta E_j\big)_-=\max\big\{\!\!-\Delta E_j,0\big\}$ is the negative part of $\Delta E_j.$ Hence, using Stokes' theorem again,
\be
\int_\Ar \big|\Delta E_j\big|=\int_\Ar \!\!\Delta E_j +2\int_\Ar \!\!\big(\Delta E_j\big)_-
=\int_{ \partial \Ar} \!\!\nabla E_j\cdot\bN+2\int_\Ar \!\!\big(\Delta E_j\big)_-
\ee
From \eqref{reg2}, \eqref{reg3} (and since $r,R$ are permissible)
\be\lim_{j\to\infty} \int_{ \partial \Ar} \!\!\nabla E_j\cdot\bN=\!\!\!\mathop\int\limits_{ \partial \Ar\setminus\cut_x}\!\! \!\!\!\nabla E\cdot\bN=A_x(R)\frac{R}{2t} E(R,t)-A_x(r)\frac {r}{2t} E(r,t)\le \C R^n t^{-{\frac n 2}-1}.\ee
Combining the above estimates we obtain
\be\ba
&\limsup_j\int_\Ar \!\!(H-\smi E_j)\Delta(H-\smi E_j)\cr&\le \limsup_j\Big( \Cr{21}\, \wt\L_x(r)\, E(r,3t)+\frac{\Cr{22}} j\Big)\int_\Ar \big(|\Delta H|+|\Delta E_j|\big)\cr&\le \C\limsup_j\Big( \wt\L_x(r)\, E(r,3t)+\frac1 j\Big)\left[\int_\Ar\!\! \Big( t^{-{\frac n 2}-1}+\frac1 j\Big)+\int_{\p\Ar}\!\! \nabla E_j\cdot\bN\right]\cr&\le
\Cl{25} R^n t^{-{\frac n 2}-1}\, \wt\L_x(r)\, E(r,3t)=\Cl{25} R^n t^{-n-1}\, \wt\L_x(r)\exp\Big(\!\!-\frac{r^2}{12t}\Big).\label{i3}
\ea\ee
To estimate the boundary integral in \eqref{fatoustokes}, using \eqref{reg2}, \eqref{gradh1}, \eqref{gradE1}
\be\ba&\lim_{j\to\infty}\int_{\p \Ar} \!\!(H-\smi E_j)\nabla (H-\smi E_j)\cdot\bN=\int_{\p \Ar\setminus\cut_x} \!\!(H-\smi E)\nabla (H-\smi E)\cdot\bN \cr&
\le\int_{\p \Ar\setminus\cut_x} \!\! |H-\smi E|\,\big(|\nabla H|+\smi|\nabla E|\big)\le \C \int_{\p \Ar}\! \wt\L_x\big(d(x,y)\big)\, E(d,3t) t^{-{\frac n 2}-\frac1 2}\cr&\le 
\C R^{n-1} t^{-{\frac n 2}-\frac1 2}  \wt\L_x(r)E(r,3t)=\C R^{n-1} t^{-n-\frac1 2}\, \wt\L_x(r)\exp\Big(\!\!-\frac{r^2}{12t}\Big).\label{i4}
\ea\ee

\medskip
Putting together \eqref{fatoustokes}, \eqref{i3}, \eqref{i4} we get, for all permissible $r,R$ with $R>r>1$, 

\be\ba &\int_\Ar \big|\nabla H-\smi\nabla E|^2\le\cr&\le\limsup_j\int_\Ar \!\!(H-\smi E_j)\Delta(H-\smi E_j)d\mu+\lim_{j\to\infty} \int_{\p \Ar} \!\!(H-\smi E_j)\nabla (H-\smi E_j)\cdot \bN\cr& \le 
\C R^{n-1} t^{-n-\frac1 2} \wt\L_x(r)\big(R t^{-\frac1 2}+1\big)\exp\Big(\!\!-\frac{r^2}{12t}\Big).\label{avest}
\ea\ee

Now note that \eqref{avest} is actually true for \emph{all} $r,R$ with $R>r>1, $ since for all such $r,R$ the integral on the left hand side of \eqref{avest} is finite (as it increases as the annulus $\Ar$ gets larger), and therefore it's separately continuous in $r $ and $R$, as well as the right hand side of \eqref{avest}.

Taking  $R=(1+\eta)r$ we obtain  \eqref{L2estH}.

 \qed

 \section{Uniform Green functions estimates}
Let us assume once again that $M$ has nonnegative Ricci curvature. Let $G_\a(x,y)$ denote the minimal, positive, Green function for $(-\Delta)^{{\frac \a 2}}$ for $\alpha$ even  and $0<\alpha<n$, i.e. a positive $C^\infty$ function off the diagonal such that $(-\Delta_y)^{{\frac \a 2}}G_\a(x,y)=\delta_x(y)$. In terms of the heat kernel we have
\be G_\a(x,y)=\frac{1}{\Gamma(\frac{\a}2)}\int_0^\infty t^{\frac\a2-1} H(x,y,t)dt\label{GH}\ee
provided the integral on the right hand side converges. Note that with $E(x,y,t)$ defined in \eqref{E} we have

\be \frac{1}{\Gamma(\frac{\a}2)}\int_0^\infty t^{\frac\a2-1} E(x,y,t)dt=c_\a d^{\a-n}(x,y),\qquad 
 c_\a=\frac{2^{-\a}\Gamma\left(\frac{n-\a}2\right)}{\pi^{n/2}\Gamma(\frac{\a}2)}\label{GE}\ee
which reduces to the classical Riesz kernel  on the standard  $\R^n$.

Arguing as in \cite{ly}, Thm 5.2, using the rough bound \eqref{hrough},   the integral in \eqref{GH} is  convergent if and only if 
\be \int_1^\infty \frac{t^{\frac\a2-1}}{V_x(\sqrt t)}dt<\infty\label{GH1}\ee
(split the integral in \eqref{GH}  at $r^2=d^2(x,y)$, make the change $s=r^4/t\ge r^2$ and use the Bishop comparison theorem in the form $V_x(r^2/\sqrt s)/V_x(\sqrt s)\ge (r^2/\sqrt s)^n/(\sqrt s)^n$).

It's easy to see that if \eqref{GH1} holds for one $x$, then it holds for all $x.$

When $\alpha=2$ condition \eqref{GH1} is equivalent to the non-parabolicity of a complete manifold with non negative ricci curvature $M$, i.e. the fact that $M$ admits a positive Green's function (see \cite{lt}, \cite{va}, \cite{ly})

\medskip
The argument that leads to \eqref{GH1} gives the rough   bound
\be \C \int_{d^2(x,y)}^\infty \frac{t^{\frac\a2-1}}{V_x(\sqrt t)}dt\le G_\a(x,y)\le  \C \int_{d^2(x,y)}^\infty \frac{t^{\frac\a2-1}}{V_x(\sqrt t)}dt.\label{GH2}\ee

\bigskip Define, for $0<\a<n$,
\be \Gt_\a(x,y)=\frac{1}{\Gamma(\frac{\a+1}2)}\int_0^\infty t^{\frac{\a-1}2}\,\nabla_y H(x,y,t)dt\label{GTH}\ee
In view of the gradient bound \eqref{gradh} it's easy to see that the above integral converges under  \eqref{GH1}, and
\be |\Gt_\a(x,y)|\le \C \int_{d^2(x,y)}^\infty \frac{t^{\frac\a2-1}}{V_x(\sqrt t)}dt.\label{GTH2}\ee
 If in addition $\a+1<n$ then we  have $\Gt_\a(x,y)=\nabla_y G_{\a+1}(x,y)$.

\smallskip We have, for $0<\a<n,$
\be \frac{1}{\Gamma(\frac{\a+1}2)}\int_0^\infty t^{\frac{\a-1}2} \nabla_y E(x,y,t)dt=-\ct_\a d^{\a-n}(x,y)\nabla_y d(x,y),\qquad 
  \ct_\a=\frac{2^{-\a}\Gamma\left(\frac{n-\a+1}2\right)}{\pi^{n/2}\Gamma(\frac{\a+1}2)}.\label{GTE}\ee
Note that
\be \ct_\a=\bc (n-\a-1)c_{\a+1}&{\text{ if }} \a+1<n\\ \dfrac{c_{\a-1}}{\a-1} &{\text { if }} 1<\a<n.\ec\label{GTE1}\ee

\bigskip

Assuming that  $M$ has Euclidean volume growth then  \eqref{GH1} is verified, and \eqref{GH2}, \eqref{GTH2} give immediately the following global rough bounds:
 \be \Cl{ggl} d(x,y)^{\a-n}\le G_\a(x,y)\le \Cl{ggu} d^{\a-n}(x,y),\quad |\Gt_\a(x,y)|\le 
\C d^{\a-n}(x,y),\quad\forall x,y\in M\label{globalG}\ee
\smallskip
In particular, for $\alpha=1$ \eqref{globalG} gives
\be|\nabla_y G_2(x,y)|\le \C d^{1-n}(x,y),\qquad \forall x,y\in M.\label{globalgradg}\ee

Finally, assuming \eqref{GH1}, we have the representation formulas for $u\in C_0^\infty(M)$
   \be u(x)=\int_M G_\a(x,y)(-\Delta)^{{\frac \a 2}} u(y)d\mu(y),\qquad   \a {\hbox { even}} \label{urep1}\ee
    \be u(x)=\int_M \lan \Gt_\a(x,y), \nabla(-\Delta)^{{\frac {\a-1} 2}} u(y)\ran d\mu(y),\qquad \a {\hbox { odd}} \label{urep2}\ee
which can be easily obtained using integration by parts, the heat equation, $E(x,\cdot, 0)=\delta_x(\cdot)$, and the fact that $t^{\frac\a2-1}V_x(\sqrt t)^{-1}\to0$ as $t\to\infty$ (a consequence of \eqref{GH1}).

\subsection{Uniform asymptotic bounds: small distances}\label{greensmall}

\bigskip

\smallskip
\begin{prop}\label{G1} Let  $M$ have nonnegative Ricci curvature and 
\be \sup_{x\in M} \int_1^\infty \frac{t^{\frac\a2-1}}{V_x(\sqrt t)}dt<\infty\label{GH3}.\ee
\smallskip \ni a) If $M$ has $1^{st}$ order bounded geometry, then there exists $\Cl{G1}>0$ such that for $d(x,y)\le 1$
\be\left|G_\a(x,y)-c_\a\,d^{\a-n}(x,y)\right|\le \Cr{G1} \, d^{\a-n+1}(x,y) \label{gsharp}.\ee

  \smallskip
  \ni b)  If $M$ has $2^{nd}$ order bounded geometry, then there exists $\Cl{Gt1}>0$ such that for $d(x,y)\le 1$ and $y\in M_x$
\be\left|\Gt_\a(x,y)+\ct_\a\,d^{\a-n}(x,y)\nabla_y d(x,y)\right|\le \Cr{Gt1} \, d^{\a-n+1}(x,y)\label{gradgsharp} .\ee
  
  \end{prop}
  
  \medskip
  \begin{rk}\label{rem-x1} When $M$ has EVG, then \eqref{GH3} holds.\end{rk}
  
  \bigskip
  
  \ni{\bf Proof.}  Assuming $1^{st}$ order bounded geometry, from  \eqref{GH}, \eqref{GE},  \eqref{hsmall}, we get 

\be \ba &\left|G_\a(x,y)-c_\a\,d^{\a-n}(x,y)\right|\le\frac{1}{\Gamma(\frac{\a}2)} \int_0^{1}t^{\frac\a2-1} \big|H(x,y,t)-E(x,y,t)\big|dt\cr&\hskip2em   +\frac{1}{\Gamma(\frac{\a}2)}\int_{1}^\infty t^{\frac\a2-1} H(x,y,t)dt+\frac{1}{\Gamma(\frac{\a}2)}\int_{1}^\infty t^{\frac\a2-1} E(x,y,t)dt\cr&
\le \frac{\Cr{4}}{\Gamma(\frac{\a}2)}  \int_0^{1}t^{-\frac n 2+\frac12+\frac\a2 -1} \exp\left(-\frac{d^2(x,y)}{24t}\right)dt+\frac{\Cr2}{{\Gamma(\frac{\a}2)} }\int_{1}^\infty\frac{t^{\frac\a2-1}}{V_x(\sqrt t)}dt+\frac{(4\pi)^{-\frac n2}}{{\Gamma(\frac{\a}2)} }\int_{1}^\infty t^{\frac\a2-\frac n 2 -1}dt
\ea\ee
and \eqref{gsharp} follows from \eqref{GH3} and 
 \be \int_0^\infty t^{-\la-1} e^{- \rho/ t}dt=\rho^{-\la}\Gamma(\la)\label{mellin}\ee
  
\medskip
  
\ni   The proof of \eqref{gradgsharp} under $2^{nd}$ order bounded geometry  is similar, using \eqref{GTH}, \eqref{GTE},  and \eqref{gradhsmall}.
   \qed

   \subsection{Uniform asymptotic bounds: large distances from a fixed point}\label{greenlarge}

\bigskip
\begin{prop}\label {G2}
 If $M$ has EVG, then  there exist $\Cl{G2}$, $\Cl{G3}$, depending only on  $n,\sm$,  such that for $x,y\in M$,  $d(x,y)\ge1$
\be\left|G_\a(x,y)-\frac{c_\a}{\sm}\,d^{\a-n}(x,y)\right|\le \Cr{G2} \,\wt\L_x\big(d(x,y)\big)d^{\a-n}(x,y).\label{gsharp1}\ee
and for any  $\eta>0$ and any  $r\ge1$ 
 \be \bigg(\fint_{A_{r,(1+\eta)r}} \Big|\Gt_\a(x,y)+\frac{\ct_\a}{\sm}d^{\a-n}(x,y)\nabla_y d(x,y)\Big|^2 dy\bigg)^{\frac1 2}\!\!\le \Cr{G3}\Big(1+\frac{1}{\sqrt\eta}\Big)
 \wt\L_x^{\frac12}(r)r^{\a-n}.\label{gradgsharp1}
 \ee
 \end{prop}
 
 \bigskip\begin{rk}\label{rem2} Assuming that  $M$ has EVG, then the global rough bounds in \eqref{globalG} imply  that \eqref{gsharp}, \eqref{gradgsharp},  \eqref{gsharp1}, \eqref{gradgsharp1}  hold for all $y\in M$.

 \end{rk}
 
 \medskip
\begin{rk}\label {rem3} When $\alpha=1$ and $n\ge3$, and under EVG
 estimate \eqref{gradgsharp1} becomes
\be \bigg(\fint_{A_{r,(1+\eta)r}} \Big|\nabla_y G_2(x,y)- \frac{(2-n)c_2}\sm d^{1-n}(x,y)\nabla_y d(x,y)\Big|^2 dy\bigg)^{\frac1 2}\!\!\le \Cr{G3}\Big(1+\frac{1}{\sqrt\eta}\Big)
 \wt\L_x^{\frac12}(r) r^{1-n}\label{gradgsharp2}
 \ee
 (in our notation $G_2(x,y)$ is the Green function for $-\Delta$),  a version  of  \cite{cm1}, eq. (4.7), with a more explicit bound in terms of the volume ratio remainder   
  (we note that in \cite{cm1} they use a different normalization for the Green function, and also that in (4.7) there should be $s^{2-2n}$).
  Compare also \eqref{gradb1} with \cite{cm1}, (3.38).
 
Without further restrictions on the geometry  one cannot expect a pointwise bound for $\wt G_\alpha$, in particular a bound of type 
\be\left|\nabla_y G_\a(x,y)-\frac{c_\a}{\sm}\,\nabla_y d^{\a-n}(x,y)\right|\le \C \,\wt\L_x\big(d(x,y)\big)d^{\a-n-1}(x,y),\ee
 which would imply \eqref{gradgsharp1}. When $\alpha=2$,  such an estimate was proved in a weaker form  in \cite{cm1}, Prop. 4.1,  under quadratic curvature decay, and we believe that it should still hold for general $\alpha$ under such assumption.
 \end{rk}
 
 \bigskip
 \ni{\bf Proof.} Estimate \eqref{gsharp1} follows at once from  \eqref{GH}, \eqref{GE}, \eqref {hlarge}, \eqref{mellin}. Likewise,   \eqref{GTH}, \eqref{GTE} imply 
  \be\ba\bigg(\fint_{A_{r,(1+\eta)r}} &\Big|\Gt_\a(x,y)+ \frac{\ct_\a}{\sm}d^{\a-n}(x,y)\nabla_y d(x,y)\Big|^2 dy\bigg)^{\frac1 2}\cr&
  \le \frac{1}{\Gamma(\frac{\a+1}2)}\int_0^\infty t^{\frac{\a-1}2}\bigg(\fint_{A_{r,(1+\eta)r}} \Big|\nabla_y H(x,y,t)- \frac{1}{\sm}\nabla_y  E(x,y,t)\Big|^2 dy\bigg)^{\frac1 2} dt \ea\ee
and \eqref{gradgsharp1} follows from   \eqref{L2estH}  and \eqref{mellin}.

\qed

\subsection{Uniform global rearrangement bounds: Talenti's type estimates}\label{Tal}

\bigskip

In this section we will derive sharp global estimates for the symmetric decreasing rearrangement of $G_\a$ in terms of the ones for the Riesz kernel $N_\a(\xi)=c_\a|\xi|^{\a-n}$ on $\R^n$, as well as similar integral bounds for $|\nabla G_\a|^q$. These results are in the spirit of the classical Talenti's comparison estimates  \cite{t1}, but in the context of manifolds with Euclidean volume growth. We will make use of recent results in this direction  due to Chen-Li \cite{cl}.

Let us recall the definition of decreasing rearrangement on a general measure space $(M,\mu)$. Given a measurable function $f\ :\ M\rightarrow \ [-\infty,\infty]$ its distribution function is defined as
\be \lambda_f(s)=\mu(\{x\in M:\ |f(x)|>s\}),\ \ \ s\geq 0.\nt\ee  
Assuming that the distribution function of $f$ is finite for all $s>0 $, the 
 decreasing rearrangement of $f$ is the function $f^*:[0,\infty)\to[0,\infty]$ defined as
\be f^*(t)=\inf\{s\geq0:\ \lambda_f(s)\leq t\},\ \ t\ge0,\label{star}\ee
where clearly $ f^*(t)=0$ if $\mu(M)\le t<\infty$.
The Schwarz symmetric decreasing rearrangement of a measurable function $f$  is the radially symmetric decreasing function $f^\#:\R^n\to\R$ defined as 
\be f^\#(\xi)=f^*(B_n |\xi|^n),\qquad\xi\in\R^n.\ee

Given a measurable function $k(x,y)$ on  $M\times  M,$ we will use the following notation for the rearrangements  of the $x-$sections and $y$-sections of $k$:
\be k^*(x,t)=\big(k(x,\cdot)\big)^*(t),\qquad  k^*(t,y)=\big(k(\cdot,y)\big)^*(t),\qquad x,y\in M,\; t\ge0.\ee


Note also that 
\be \lambda_f(s)\le A s^{-\sgam},\quad \forall s>0\;\;\iff\;\; f^*(t)\le A^{\frac1\sgam}t^{-\frac1\sgam},\quad \forall t>0\ee

and if either inequality holds a.e. then it also holds everywhere.
\bigskip
\begin{theorem}\label{talenti} Let $M$ have EVG, $\ric\ge0$, {\rm{dim}}$M=n\ge3$, and denote 
\be \lambda_\a(x,s)=\mu(\{y: G_\a(x,y)>s\}),\qquad x\in M.\label{df}\ee
Then, for any $\alpha$ even, $2\le\alpha<n$, and  all $x\in M$, we have 
\be G_\a^*(x,t)\le \sm^{-\frac\a n} N_\a^*(t)=\sm^{-\frac\a n} B_n^{\frac{n-\a}n} c_\a \;t^{-\frac{n-\a}n},\qquad  t>0,\label{Tb}\ee
and for  $0<q\le2$,  $0<s_1<s_2$, and all $x\in M$, we have 
\be\mathop\int\limits_{s_1<G_{\a}(x,\cdot)\le s_2}\hskip-1.5em  |\nabla_y G_{\alpha}(x,y)|^q dy\le \sm^{-\frac {q(\a-1)}n } \hskip-2.5em \mathop\int\limits_{\lambda_{\a}(x,s_2)<B_n|\xi|^n\le\lambda_{\a}(x,s_1)}\hskip-2em \big |\nabla N_\a(\xi)\big|^qd\xi.
\label{Tc}\ee
In particular,  for all $x\in M$, $0<s_1<s_2$, $\alpha$ even, and  $2\le\alpha<n$, we have
\be\hskip-1em\mathop\int\limits_{s_1<G_{\a}(x,\cdot)\le s_2} \hskip-1em |G_{\alpha}(x,y)|^\nna dy\le\sm^{-\frac\a{n-\a}}\gamma_{n,\a}^{-1}\log\frac{\lambda_{\a}(x,s_1)}{\lambda_{\a}(x,s_2)}
\label{Td}\ee
whereas for  $\alpha$  odd and  $1\le\a\le \frac n2$ 
\be\hskip-.51em \mathop\int\limits_{s_1<G_{\a+1}(x,\cdot)\le s_2}\hskip-1.5em  |\nabla_y G_{\alpha+1}(x,y)|^\nna dy\le\sm^{-\frac\a{n-\a}}\gamma_{n,\a}^{-1}\log\frac{\lambda_{\a+1}(x,s_1)}{\lambda_{\a+1}(x,s_2)}.
\label{Te}
\ee

\end{theorem}
\bigskip
\begin{rk}\label{r1} As it will be apparent from the proof below, we will also obtain the sharp bounds \eqref{Tb}-\eqref{Te} for the Green functions  of the Poisson equation  with homogeneous Navier boundary condition on any smooth bounded domain $\Omega\subseteq M$. Specifically, given 
 \be\bc (-\Delta)^{\frac\a2}u=f \text{ on } \Omega\cr  (-\Delta)^j u=0 \text{ on } \partial \Omega,\; 0\le j\le \lceil (\a-1)/2.\rceil\ec\label{navier}\ee
and the corresponding symmetrized problem 
\be\bc (-\Delta)^{\frac\a2}v=f^\# \text{ on } \Omega^\#\cr  (-\Delta)^j v=0 \text{ on } \partial \Omega^\#,\; 0\le j\le \lceil (\a-1)/2.\rceil\ec\label{navier1}\ee
where $\Omega^\#$ denotes the ball of center 0  in $\R^n$ and with volume $\mu(\Omega)$, if 
 $G_\a^{\Omega}(x,y)$ denotes the Green function for \eqref{navier}, and $N_\a^{\Omega^\#}(\xi)$ the  Green function for
 \eqref{navier1} with pole $0$, then \eqref{Tb}-\eqref{Te} hold with $G_\a^\Omega$ in place of $G_\a$ and $N_\a^{\Omega^\#}$ in place of $N_\a$, where $\lambda_\a$ is of course defined with respect to $\Omega.$

When $M=\R^n$ ($\sm=1$) bounds \eqref{Tb} for the Dirichlet Green function $G_2^\Omega$ appear in Bandle \cite{ban}, whereas for the biharmonic Green function ($\alpha=4$)  they can be deduced from results in \cite{fk}. The integral bounds \eqref{Tc}, \eqref{Td}, \eqref{Te} do not seem to appear in the literature, even on $\R^n$. Note that letting $s_2\to\infty$ and $s_1\to0$ in \eqref{Tc}, allows us to replace the domains  of integration with $\Omega$ and $\Omega^\#$ respectively. Such an estimate on $\R^n$ was given by Weinberger \cite{w} for $\alpha=2$ and for $q<\frac n{n-1}$ (see also \cite{ban}, Thm 2.5).

\bigskip
\ni{\bf Proof.}   We will begin to prove \eqref{Tb}, \eqref{Tc}  for the Dirichlet Green function $G_2^\Omega$ of a smooth bounded domain $\Omega$ in $M$. Let us  recall some classical results first. If $\Omega_1\subseteq \Omega_2$ then $G_2^{\Omega_1}(x,y)\le G_2^{\Omega_2}(x,y)$ for $x,y\in\Omega_1$, $x\neq y$,  by the maximum principle and the local asymptotic expansion of $G_2^\Omega$ (see e.g. \cite{sy}, p. 82).  Secondly, under our assumptions $M$ is non parabolic, and if $\{\Omega_k\}$ is an increasing sequence of smooth, bounded,  open sets  which exhaust $M$, then for any fixed $x,y\in M$, ($x\neq y$) we have $x,y\in \Omega_k$ for $k\ge k_0$, and $G_2^{\Omega_k}(x,y)\uparrow G_2(x,y)$, where $G_2$ is the unique positive, minimal Green function defined earlier. 

Recall now  a recent result by Chen-Li \cite{cl}, which gives a version of Talenti's comparison theorem on a noncompact Riemannian manifold $M$, with EVG and $\ric\ge0$. If $\Omega$ is a smooth, bounded open set in $M$, and $f\in L^2(\Omega)$ then let $u\in W_0^{1,2}(\Omega)$ be the weak solution of 
\be\bc -\Delta u=f & {\text {in }} \Omega\\ u=0  & {\text {on }} \partial \Omega\ec\label{Dir1}\ee
and let $v\in W_0^{1,2}(\Omega^\#)$ be the weak solution of
\be\bc -\Delta v=f^\# & {\text {in }} \Omega^\#\\ v=0  & {\text {on }} \partial \Omega^\#.\ec\label{Dir2}\ee

Then, Chen-Li proved that \be u^*(t)\le \sm^{-\frac2n} v^*(t),\qquad t>0.\label{CL} \ee

The proof is based on Talenti's original argument, in the version given by Kesavan \cite{k}, combined with the sharp isoperimetric inequality on manifolds with EVG. In particular, if 
\be \lambda(s):=\mu(\{|u|>s\}),\qquad F(t)=\int_0^t f^*(\zeta)d\zeta\ee then
\be1\le \sm^{-\frac2n} (nB_n^{\frac1n})^{-2} \lambda(s)^{\frac2n-2} F\big(\lambda(s)\big)(-\lambda'(s)),\qquad {\hbox{ a.e. }} s>0.\label{T1}\ee

Arguing as in \cite{k}, p. 51, it's easy   to see that for $0<q\le 2$
\be {-\frac{d}{ds}}\int_{|u|>s} |\nabla u|^q dx\le  F\big(\lambda(s)\big)^{\frac q2}(-\lambda'(s))^{\frac{2-q}2},\qquad {\hbox{ a.e. }} s>0. \label{T2}\ee

Multiplying the  last estimate by \eqref{T1} raised to the power $q/2$, integrating  in $s\in[s_1,s_2]$, and changing variables yields
\be\ba \int_{s_1<|u|\le s_2}&|\nabla u|^q dx\le \sm^{-\frac q n}\int_{s_1}^{s_2}\big((n B_n^{\frac 1 n})^{-1}\lambda(s)^{\frac1n-1} F(\lambda(s))\big)^q (-\lambda'(s)) ds\cr&\le \sm^{-\frac q n}\int_{\lambda(s_2)}^{\lambda(s_1)}\big((n B_n^{\frac 1 n})^{-1}t^{\frac1n-1} F(t)\big)^q dt =\sm^{-\frac q n}\hskip-1em\mathop \int\limits_{\lambda(s_2)<B_n|\xi|^n\le \lambda(s_1)}\hskip-1.5em|\nabla v(\xi)|^q d\xi,\label{T3}\ea\ee
where for the last identity we used the explicit formula
\be v(\xi)=n^{-2}B_n^{-\frac2n}\int_{B_n|\xi|^n}^{|\Omega^\#|} \tau^{\frac2n-2}F(\tau)d\tau.\ee

With this in mind, we can readily obtain \eqref{Tb},\eqref{Tc} for the Dirichlet Green function $G_2^\Omega$ of  a smooth domain $\Omega$ inside $M$.

Fix $x\in \Omega $, and let $f_k\in C_c^\infty(M)$ be such that $f\ge0$,  $\int f_k=1$ and $\supp f_k\subseteq B(x,{1\over k})\subseteq \Omega$, for $k\ge k_0$. Then $f_k\to \delta_x$ in the sense of distributions. Let $u_k$ be the solution of  the Dirichlet problem \eqref{Dir1} with  $f=f_k$, and let $v_k$ be the solution of the corresponding symmetrized problem \eqref{Dir2}. The function $u_k$ is given by the formula
\be u_k(y)=\int_\Omega G_2^\Omega(y,z)f_k(z)dz,\qquad y\in\Omega\label{NP}\ee
and it's smooth. Fix $y\in \Omega$, $y\neq x$ and assume $\frac1{k_0}<\frac12 d(y,x)$. If  $\varphi$ is a nonnegative bump function and supported on $B(x,\frac1{k_0})$ and equal 1 on $B(x,\frac1{k_0+1})$, then for $k>k_0+1$
\be u_k(y)=\int_{\Omega} G_2^\Omega(y,z)\varphi(z)f_k(z)dz\to u(y):=G_2^\Omega(y,x)\label{uk1}\ee
and moreover, since $G_2^\Omega(y,z)\le G_2(y,z)\le \Cr{ggu} d(y,z)^{2-n}$ (from \eqref{globalG}), 
\be0\le u_k(y)\le \int_{B(x,\frac1{k_0})} G_2^\Omega(y,z)\varphi(z)f_k(z)dz\le \Cr{ggu} \int_{B(x,\frac1{k_0})} d(y,z)^{2-n}f_k(z)dz\le 2^{n-2}\Cr{ggu} d(y,x)^{2-n}.\label {uk2}\ee

Regarding $v_k$,  we have the explicit formula

\be v_k^*(t)=n^{-2}B_n^{-\frac2n}\int_t^{|\Omega^\#|}\tau^{\frac2n-2}F_k(\tau)d\tau ,\qquad F_k(\tau)=\int_0^\tau f_k^*\ee
valid for $0<t\le|\Omega^\#|=\mu(\Omega)$, where the rearrangement is made with respect to $\Omega$. Since $\int f_k=1$ we then have $F_k\le 1$, and in fact for any $t\in (0,|\Omega^\#|]$, we have $F_k(t)=1$ for $k$ large enough (so that $\mu(|\supp f_k|)<t$).

Hence,
\be v_k^*(t)\le \frac1{n(n-2)B_n^{\frac2n}}\big(t^{-\frac {n-2}n}-|\Omega^\#|^{-\frac{n-2}n}\big):=(N_2^{\Omega^\#})^*(t)\ee
where, recall, $N_2^{\Omega^\#}$ is the  Dirichlet Green function for the ball of volume $|\Omega^\#|=\mu(\Omega)$, and  with pole $0$.

By \eqref{CL} we have, for $0<t\le |\Omega^\#|$
\be u_k^*(t)\le \sm^{-\frac2n} v_k^*(t)\le\sm^{-\frac2n}(N_2^{\Omega^\#})^*(t) \le \frac{\sm^{-\frac2n}}{n(n-2)B_n^{\frac2n}}\;t^{-\frac {n-2}n}. \ee
and hence
\be \lambda_{u_k}(s)\le \bigg( \frac{\sm^{-\frac2n}}{n(n-2)B_n^{\frac2n}}\bigg)^{\frac n{n-2}}s^{-\frac n{n-2}},\qquad s>0.\ee
Now, since $\lambda_u$ is discontinuous at countably many points, then $\mu(\{y: |u(y)|=s\})=0$ for a.e. $s>0$. For any such $s$ we have that    $\chi_{\{|u_k|>s\}}\to \chi_{\{|u|>s\}}$,  $\mu$-a.e., and using the Dominated Convergence Theorem and the fact that $\mu(\Omega)<\infty$, we can conclude that $\lambda_{u_k}(s)\to\lambda_{u}(s)$. Note that the convergence is actually true for all $s>0$, in the case of the Green function,  since by classical  results (see for ex. \cite{che}, Thm. 2.2)  we have 
$\mu(\{z\in\Omega:G_2^\Omega(y,z)=s\})=0$ for all $s>0$ (see also Remark \ref{cheng}). Thus we can conclude that 
\be \lambda_{u_k}(s)\le \bigg( \frac{\sm^{-\frac2n}}{n(n-2)B_n^{\frac2n}}\bigg)^{\frac n{n-2}}s^{-\frac n{n-2}}\ee
for a.e. $s>0$ and hence for all $s>0.$ The above inequality implies then \eqref{Tb}, when $\alpha=2.$, and for the Green function $G_2^\Omega$.

The proof of \eqref{Tc} for $\alpha=2$  follows similarly from \eqref{T3} applied to $u_k, v_k$. Indeed, using \eqref{NP} we can see easily  that $|\nabla_y u_k(y)|^2\to|\nabla_y u(y)|^2=|\nabla_y G_2^\Omega(y,x)|^2$ for $y\in\Omega$ and Fatou's lemma combined with \eqref{T3} gives \eqref{Tc}.

The proof of \eqref{Tb}, \eqref{Tc} for general $\alpha$ even, is done by iteration.\, using the relations
\be G_\a^\Omega(x,y)=\int_M G_2^\Omega(x,z)G_{\a-2}^\Omega(z,y)dz,\qquad -\Delta_y G_\a^{\Omega}(x,\cdot)=G_{\alpha-2}^\Omega(x,\cdot)\label{iter}\ee
where $G_\a^\Omega$ denotes the  Green function for the problem $(-\Delta)^{\frac\a2}u=f$ on $\Omega$, with homogeneous Navier boundary condition $(-\Delta)^j u=0$ on $\partial \Omega$, for $0\le j\le \lceil (\a-1)/2\rceil$.

Indeed,  fix $x\in M$ and let $u_{k,\a},v_{k,\a}$ be the solutions to the problems
\be\bc -\Delta u_{k,\a}=u_{k,\a-2} & {\text {in }} \Omega\\ u_{k,\a}=0  & {\text {on }} \partial \Omega\ec\qquad \bc -\Delta v_{k,\a}=u_{k,\a-2}^\# & {\text {in }} \Omega^\#\\ v_{k,\a}=0  & {\text {on }} \partial \Omega^\#\ec\ee
where $u_{k,2}=f_k$.

Assuming that for each $y\in \Omega$
\be u_{k,\a-2}(y)\to G_{\a-2}^\Omega(y,x),\qquad   \text {as } k\to\infty\label{estuk1}\ee
\be 0\le u_{k,\a-2}(y)\le C_\a d(x,y)^{\a-2-n},\qquad  \text{for } k>2k_0\label {estuk2}\ee
and 
\be u_{k,\a-2}^*(t)\le \sm^{-\frac{\a-2}n} B_n^{\frac{n-\a-2}n} c_{\a-2} t^{-\frac{n-\a-2}n},\qquad t>0, \label{estuk3}\ee
where $C_\a$ is a constant independent of $x,y,\Omega$, 
by induction we can show that  see that \eqref{estuk1}, \eqref{estuk2}, \eqref{estuk3}, and hence \eqref{Tb},  also hold for $u_{k,\a}$.

Indeed for $\alpha= 4$, estimate \eqref{estuk2} is is just \eqref{uk2}, whereas for $\alpha>4$ we have

 \be\ba0\le u_{k,\a}(y)&= \int_\Omega G_2^\Omega(y,z) u_{k,\a-2}(z)dz\le C_\a\int_{\Omega}G_2^\Omega(y,z) d(z,x)^{\a-2-n}dz\cr &\le C_{\a}\Cr{ggl}\int_M G_2(y,z)G_{\a-2}(z,x)dz=C_\a\Cr{ggl} G_\a(y,x)\le C_{\a+2} d(y,x)^{\a-n}.\ea\ee
where we used that $G_2^\Omega\le G_2$, and the global bounds \eqref{globalG}. The proof of \eqref{estuk1} for $u_{k,\a}$  follows by induction using the Dominated Convergence Theorem. The proof of \eqref{estuk3} for $u_{k,\a}^*$  follows also by induction, since by Talenti's estimate and
\be\ba u_{k,\a}^*(t)&\le\sm^{-\frac2n} v_{k,\a}^*(t)=\sm^{-\frac2n} n^{-2}B_n^{-\frac2n}\int_t^{|\Omega^\#|}\tau^{\frac2n-2}\int_0^\tau u_{k,\a-2}^*(\zeta)d\zeta\cr&\le\sm^{-\frac2n} n^{-2}B_n^{-\frac2n}\int_t^{|\Omega^\#|}\tau^{\frac2n-2}\int_0^\tau 
\sm^{-\frac{\a-2}n} B_n^{\frac{n-\a+2}n}c_{\a-2}\zeta^{-\frac{n-\a+2}n}d\zeta\cr&=
\sm^{-\frac\a n}B_n^{\frac{n-\a}n}\frac{c_{\a-2}}{(\a-2)(n-\a)} \big(t^{-\frac{n-\a}n}-|\Omega^\#|^{-\frac{n-\a}n}\big) 
\ea\ee
and the result follows from the identity $c_\a{(\a-2)(n-\a)}={c_{\a-2}}$.

To prove  \eqref{Tc} for $2<\alpha\le \frac n2$, use the relation
\be u_{k,\a}(y)=\int_{\Omega} G_2^\Omega(y,z)u_{k,\a-2}(z)dz.\ee
Note that from the classical gradient bound for positive harmonic functions (see e.g. \cite{li}, Theorem 6.1) we get that for all $y,z\in \Omega$, $y\neq z$
\be |\nabla_y G_2^\Omega(y,z)|\le C(n) G_2^\Omega (y,z) d(y,z)^{-1}\le C(n) G_2(y,z)d(y,z)^{-1}\le C(n)\Cr{ggu}  d(y,z)^{1-n},\label{ng1}\ee
for some $C(n)$ depending only on $n$. From here, using the Dominated Convergence Theorem
one sees easily that $|\nabla_y u_{k,\a}(y)|^2\to|\nabla_y G_\a^{\Omega}(y,x)|^2$ for all $y\in \Omega$, and Fatou's 
lemma combined with \eqref{T3} applied with $u=G_\a^\Omega(x,y)$ gives \eqref{Tc}.

The proof of \eqref{Td} and\eqref{Te}, for the Dirichlet Green functions,  is a simple computation based on \eqref{Tb}, \eqref{Tc}; we leave the details  to the reader.

To prove the theorem for the global Green functions, consider an increasing sequence of smooth bounded sets $\{\Omega_k\}$ such that $\Omega_k\uparrow M$, so that the corresponding Dirichlet Green functions satisfy $G_2^{\Omega_k}\uparrow G_2$, and as a consequence,  for any $\alpha $ even $G_\a^{\Omega_k}\uparrow G_\a$. From this relation it's clear that for any $x\in M$ and $t>0$  we have $(G_\a^{\Omega_k})^*(x,t)\uparrow (G_\a)^*(x,t)$ so that \eqref{Tb} holds.

Next, note that if $\lambda_\a^k(x,s)=\mu(\{y\in\Omega_k: G_\a^{\Omega_k}(x,y)>s\})$ then $\lambda_\a^k(x,s)\uparrow \lambda_\a(x,s)$. Also,  
Using the gradient bound applied to the positive harmonic function $G_2(x,\cdot)-G_2^{\Omega^k}(x,\cdot)$ we  obtain, for fixed $x,y\in M$, $x\neq y$,
\be\big|\nabla_y\big(G_2(x,y)-G_2^{\Omega^k}(x,y)\big)\big|\le C(n) \big(G_2(x,y)-G_2^{\Omega^k}(x,y)\big) d(x,y)^{-1}\to0,\qquad k\to\infty\ee
which implies $|\nabla_y G_2^{\Omega_k}(x,y)|\to |\nabla_y G_2(x,y)|$, as $k\to+\infty$. Using that 
\be G_\a^{\Omega_k}(x,y)=\int_{\Omega_k} G_{\a-2}^{\Omega_k}(x,z)G_2^{\Omega_k}(z,y)dz\ee
the bounds $G_{\a-2}^{\Omega_k}(x,z)\le \Cr{ggu} d(x,z)^{\a-2-n}$, estimate \eqref{ng1}, and the Dominated Convergence Theorem, gives \eqref{Tc}, \eqref{Td}, \eqref{Te}.

\qed

\medskip
\begin{rk}\label{r2} As observed in Remark \ref{r1}, the previous yields \eqref{Tb}-\eqref{Te} for the Green functions of smooth bounded domains. One just has to be a bit more careful and assume, in the induction process, that \eqref{estuk3} is valid in the stronger  form $u_{k,\a-2}^*(t)\le \sm^{-\frac{\a-2}n}(N_{\a-2}^{\Omega^\#})^*(t)$, and using that
\be \big(N_{\a}^{\Omega^\#}\big)^*(t)=n^{_2}B_n^{-\frac2n}\int_{t}^{|\Omega^\#|}\tau^{\frac2n-2} \int_0^\tau  \big(N_{\a-2}^{\Omega^\#}\big)^*(\zeta)d\zeta.\ee

\section{Proofs of Theorems \ref{main1} and Theorem \ref{main2}: inequalities}

Let us recall the Adams inequalities on metric measure spaces proved recently in \cite{mq1}.

\medskip
Let $(M,d,\mu)$ be a metric measure space, that is a set $M$ endowed with a distance function $d$ and a Borel measure $\mu$.
 Define 
$$B(x,r)=\{y:\; d(x,y)\le r\},\qquad V_x(r)=\mu\big(B(x,r)\big).$$
and assume that  that for all $x\in M$
\be i)\; \; \forall r>0,  V_x(r)<\infty,\qquad ii)  \; r\to V_x(r) \;{\rm continuous.}\label{vol}\ee
In particular, the boundary of any ball has zero measure, and $\mu$ is nonatomic.

In \cite{mq1}  we defined  a  measurable $k:M\times M\to \R$ to be  {\it Riesz-like kernel of order $\beta>1$ and normalization constants $A_0>0$ and $A_\infty\ge0$} if  there exists $B\ge0$ such that for all $x\in M$

\be\mathop\int\limits_{r_1\le d(x,y)\le r_2} |k(x,y)|^\b d\mu(y)\le  A_0\log{\frac{V_x(r_2)}{V_x(r_1)}}+B,  \quad  0<V_x(r_1)<V_x(r_2)\le1,\label{k1}\tag{K1}\ee
 
  \be\mathop\int\limits_{r_1\le d(x,y)\le r_2} |k(x,y)|^\b d\mu(y)\le  A_\infty\log{\frac{V_x(r_2)}{V_x(r_1)}}+B,  \quad  1\le V_x(r_1)<V_x(r_2)\label{k2}\tag{K2}\ee

for all $x,y\in M$
\be |k(x,y)|\le BV_x\big(d(x,y)\big)^{-1/\b}\qquad |k(x,y)|\le B V_y\big(d(x,y)\big)^{-1/\b}\label{k3}\tag{K3}\ee

\medskip
for each $\delta>0$  there is $B_\delta>0$ such that for all $x\in M$

\be\mathop\int\limits_{ d(x,y)> R} |k(x,y)-k(x',y)|^\b d\mu(y)\le B_\delta,\quad V_x(R)\ge (1+\delta) V_x(r),\quad \forall x'\in B(x,r).\label{k4}\tag{K4}\ee

A Riesz-like potential is an integral operator
 \be Tf(x)=\int_M k(x,y)f(y)d\mu(y).\label{de3}\ee
 where $k(x,y)$ is  Riesz-like. Such operator is well defined for $f\in L^{\b'}$ and with compact support.  Here and from now on $\b'$ will denote the exponent conjugate to $\b$:
 $$\frac {1}\b+\frac 1{\b'}=1.$$

 \begin{theorem}[\cite{mq1}]\label{m1} If $k(x,y)$ is a Riesz-like kernel of order $\b>1$ and normalization  constants $A_0, A_\infty$, and if  
$p\ge 1$, 
 then there is a constant $C$ such that for any measurable function $f$ supported in a ball with 

 \be\|f\|_{\b'}\le 1\label{1b}\ee we have
 
 \be \int_M \frac{\exp_{\llceil p/\b-1\rrceil}\left(\dfrac{1} {\max\{{A_0},{A_\infty}\}}|Tf|^{\b}\right)}{1+|Tf|^{p\b/\b'}}\,d\mu\le C\|Tf\|_{p}^{p}. \label{1a}\ee

\def\k{\kappa}
Moreover, if $A_\infty>0$ and given any  $\k>0$ there exists $C$ (depending on $\k$) such that for any measurable  $f$ supported in a ball with 
\be \|f\|_{\b'}^{\b'}+\k \|Tf\|_{\b'}^{\b'}\le 1\label{1ba}\ee 
we have 
 \be \int_M \frac{\exp_{\llceil p/\b-1\rrceil}\left(\dfrac1 {A_0}|Tf|^{\b}\right)}{1+|Tf|^{p\b/\b'}}\,d\mu\le C\|Tf\|_{p}^{p}, \label{1a1}\ee
and
 \be \int_M \exp_{\llceil \b'-2\rrceil}\left(\dfrac1 {A_0}|Tf|^{\b}\right)\,d\mu\le  {C}. \label{1a2}\ee

 \end{theorem}
   
   \medskip
   
    \begin{rk}\label{rkvector} The conclusions of Theorem \ref{m1} hold if $k$ and $f$ are vector-valued. Specifically, if $k=(k_1,...,k_m)$ is defined on $M\times M$ and measurable, and if $f=(f_1,...,f_m)$ is measurable on $M$, then let $Tf(x)=\int_M k(x,y)\cdot f(y)d\mu$ where the $``\,\cdot\, ”$ denotes the standard Euclidean inner product on $\R^m.$ If $k$ satisfies \eqref{k1}-\eqref{k4}, with $|k|=(k\cdot k)^{1/2},\ |f|=(f\cdot f)^{1/2}$, $\;\|f\|_{\b'}=\big(\int_M|f|^{\b'}\big)^{1/\b'}$,   then the conclusions of Theorem {m1} hold.

    \bigskip
    Clearly, we would like to apply the above theorem  when $(M,\mu,d)$ is a Riemannian manifold satisfying the assumptions of Theorems \ref{main1},\ref{main2} (it is known that the conditions in \eqref{vol} are true). Given the representation formulas \eqref{urep1}, \eqref{urep2}, we want to take $T=T_\a$ to be  the integral operator with kernel $G_\a(x,y)$ for $\alpha$ even, and $\Gt_\a(x,y)$ for $\alpha$ odd. 
  What we will verify is that without bounded geometry  $G_\a$ and $\Gt_\a$ satisfy conditions \eqref{k1}-\eqref{k4}, with $\b=\nna$,  $A_0=A_\infty=
\sm^{-\frac{\a}{(n-\a)}}\gam_{n,\a}^{-1},$ and if bounded geometry is assumed then $A_0=\gam_{n,\a}^{-1}<A_\infty$.   This will imply the result for $u\in C_c^\infty(M)$, and by density  for $u\in H^{\a,\na}(M)$.

 \begin{rk}\label{dp} When $\alpha$ is odd the operator defined by \eqref{urep2} is not one of the  types described in Remark \ref{rkvector}, i.e. involving a Euclidean inner product. However, we can argue as in the proof of Theorem 7 in  \cite{mq1}, to reduce matters to Euclidean inner products. Namely, let us pick any point $p\in M$ and a chart in geodesic polar coordinates at $p$, and in the notation of Section \ref{background},   $\exp_p$ is a diffeomorphism from $D_p$ onto $M\setminus (\cut_p\cup \{p\})$, where we can identify $D_p$ with  $\big\{r\xi:\;0<r<c(\xi),\, \xi\in \S^{n-1}\big\}\subseteq \R^n$.  For any $\yb\in D_p$  we have $\big(g_{ij}(\yb)\big)=R_\yb^TR_\yb$ for some invertible matrix $R_\yb$, and 
\be u(x)=T_\a f(\xb)=\int_{D_p} K_\a(\xb,\yb)\cdot f(\yb)\sqrt{|g(\yb)|}\, dm(\yb),\qquad x=\exp_p \xb, \ee
where 
\be K_\a(\xb,\yb)=R_\yb\, \wt G_\a(x,y),\qquad f(\yb) =R_\yb \nabla D^{\a-1}u(y)=R_\yb D^\a u(y).\ee
Since $K_\a\cdot K_\a=|\wt G_\a|^2$,  $f\cdot f=|\nabla D^{\a-1}u|^2$, and $\cut_p$ has zero measure, we can check conditions (K1)-(K4) using the norm in the metric $g$, and then apply Theorem \ref{m1} and Remark~\ref{rkvector} to the operator $T_\a$  defined on the space $D_p$, endowed with metric $d(\xb,\yb)=d(x,y)$ and measure $d\mu(\yb)=\sqrt{|g(\yb)|}\,dm(\yb)$.
\end{rk}
\smallskip

    \smallskip
   We will now derive \eqref{k1}-\eqref{k4}  without any bounded geometry assumptions, based on the bounds of Theorem \ref{talenti}. Assume $n\ge 3$.  The global bounds \eqref{globalG} imply that for any fixed $x\in M$, and  any $r_2>r_1>0$
   \be \{y: r_1\le d(x,y)< r_2\}\subseteq \{y:\Cr{ggl} r_2^{\a-n}< G_\a(x,y)\le \Cr{ggu} r_1^{\a-n}\}\ee
and for the distribution function defined in \eqref{df}
\be V_x\big(\C s^{-\frac1{n-\a}}\big)\le \lambda_\a(x,s)\le  V_x\big(\C s^{-\frac1{n-\a}}\big),\qquad s>0.\ee
The above estimates combined with \eqref{Tb}, \eqref{Tc} and  \eqref{Evg3},  instantly imply that there is $B>0$ such that for any $x\in M$
\be\mathop\int\limits_{r_1\le d(x,y)\le r_2}\hskip-1.5em G_\a(x,y)^{\nna} dy\le \sm^{-\frac\a{n-\a}}\gamma_{n,\a}^{-1}\log\frac{V_x(r_2)}{V_x(r_1)}+B,\qquad 0<r_1<r_2\ee
when $\a$ is an even integer, and the same inequality holds for $\Gt_\a(x,y)$ when $\a$ is an odd integer smaller or equal $n/2.$
 This shows that \eqref{k1}, \eqref{k2} hold with $A_0=A_\infty=\sm^{-\frac\a{n-\a}}\gamma_{n,\a}^{-1}$.

   \smallskip
 Regarding \eqref{k3}, such estimate  follows  from \eqref{Evg3} and the global bounds  \eqref{globalG}. 
   
   To prove \eqref{k4},  assume $V_x(R)\ge (1+\delta)V_x(r)$, for some $\delta>0$. Since $V_x(R)/|B_\sR^n(r)|$ is decreasing and smaller than 1, then  $V_x(R)-V_x(r)\le B_n(R^n-r^n)$, from which we deduce $R^n-r^n\ge\sm\delta r^n$, and therefore $R\ge (1+\eta_\delta)r$, with $\eta_\delta=(\sm\delta)^{\frac1n}$. By integrating in~$t$ \eqref{hol1} and \eqref{hol2} we then get the following Lipschitz estimates, for $d(x,x')\le r,\; d(x,y)\ge R$:
 \be |G_\a(x,y)-G_\a(x',y)|\le C_\delta\, d(x,x') d^{\a-n-1}(x,y)\label{holderG}\ee
 \be  |\wt G_\a(x,y)-\wt G_\a(x',y)|\le \wt C_\delta\, d(x,x') d^{\a-n-1}(x,y),\label{holderG1}\ee
 from which \eqref{k4} follows easily. This proves the inequality part  of Theorem \ref{main1}, in case $n\ge3$. 
 
 When $n=2$ we only need to consider $\alpha=1$. In this case however we can use directly the P\'olya-Szeg\H o inequality \eqref{PS} to prove \eqref{MS} for $\alpha=1$ and any $n\ge2$. Indeed, assuming $\|\nabla u\|_n\leq1$ we have $\|\sm^{\frac1n} \nabla u^\#\|_n\le 1$, so using the fact that \eqref{MS} is true on $\R^n$, 
 \be\ba \int_M &\frac{\exp_{\llceil n-2\rrceil}\left(\sm^{\frac1{n-1}}\gamma_{n,1}|u|^{\frac n{n-1}}\right)}{1+|u|^{\frac n{n-1}}}\,d\mu\cr&= \int_{\R^n} \frac{\exp_{\llceil n-2\rrceil}\left(\gamma_{n,1}|\sm^{\frac1 n}u^\#|^{\frac n{n-1}}\right)}{1+|u^\#|^{\frac n{n-1}}}\,dm\le C\|u^\#\|_n^n=C\|u\|_n^n,
 \ea\label{PSZ}\ee where $dm$ denotes the Lebesgue measure on $\R^n$.
This concludes the proof the inequalities in Theorem \ref{main1}. Now let us prove \eqref{MSR} and   \eqref{MT}.

Assuming  that  $\alpha$ even, and that $M$ has  $1^{\rm {st}}$ order bounded geometry, let us show that \eqref{k1} holds with $A_0=\gamma_{n,\a}^{-1}$.  Note that if $V_x(r)\le 1$ then $r\le \big(\sm B_n\big)^{-\frac1 n}$.
    Using the inequality $|A^p-B^p|\le p|A-B|\max\{A^{p-1}, B^{p-1}\}$, ($A,B\ge 0$, $p>1$), and using \eqref{gsharp}, \eqref{globalG},   and Remark \ref{rem2}, we obtain, for $d(x,y)\le\big(\sm B_n\big)^{-\frac1 n}$,
    
    \be \Big||G_\a(x,y)|^\nna- \big(c_\a d^{\a-n}(x,y)\big)^\nna\Big|\le \Cl{A} d^{\a-n+1}(x,y) \Big(d^{\a-n}(x,y)\Big)^{\frac{\a}{n-\a}}=\Cr{A}d^{1-n}(x,y)\ee
    
    hence 
    \be\bigg|\int_{r_1\le d(x,y)\le r_2}|G_\a(x,y)|^\nna dy-c_\a^{\nna}\int_{r_1\le d(x,y)\le r_2} d^{-n}(x,y)dy\bigg|\le \C,\qquad 0<r_1<r_2\le \big(\sm B_n\big)^{-\frac1 n}\ee
    
    From the Bishop comparison \eqref{Evg3}
    
    \be\ba \int_{r_1\le d(x,y)\le r_2} d^{-n}(x,y)dy&\le\int_{r_1\le d(x,y)\le r_2} \frac{B_n}{V_x\big(d(x,y)\big)}dy\cr&={B_n}\int_{V_x(r_1)}^{V_x(r_2)}\frac 1 t dt={B_n}\log\frac{V_x(r_2)}{V_x(r_1)}\ea\label{est2}\ee
which implies \eqref{k1} with $A_0=c^{\nna}B_n$.   Note that   the first identity in \eqref{est2} follows from the fact that 
\be \bigg(\frac 1{V_x\big(d(x,\cdot)\big)}\bigg)^*(t)=\frac 1 t\ee
and the fact that balls are level sets of $V_x(d(x,y))^{-1}$ (up to sets of zero measure on general metric measure spaces satisfying \eqref{vol}).

To show \eqref{k1} in the case $\alpha$ odd, under $2^{\rm {nd}}$ order  bounded geometry, we proceed in exactly the same manner, using \eqref{gradgsharp}, and the fact that $|\nabla_y d(x,y)|=1.$

Since we have already checked \eqref{k2},\eqref{k3},\eqref{k4}, and $A_0\le A_\infty$, the inequality part  of Theorem \ref{main2} is proved, for $n\ge 3$.

\medskip
\begin{rk}\label{LMT} When $\alpha=1$ the proof of \eqref{MT} can be achieved using the results in \cite{fmq}. Indeed,\eqref{globalgradg} and  \eqref{gradgsharp} imply that 
\be |\nabla_y G_2(x,\cdot)|^*(t)\le \gamma_n^{ -\frac{n-1}n} t^{-{\frac{n-1}n}}(1+ C t^{\frac 1n}),\qquad 0<t\le 1\ee
and
\be |\nabla_y G_2(x,\cdot)|^*(t)\le C t^{-{\frac{n-1}n}},\qquad t>0,\ee
and the constant $C$ in both inequalities is independent of $x$.
According to \cite{fmq}, Observation 3, we can conclude that the Moser-Trudinger inequality holds locally, i.e. for functions supported on sets with measure less than 1, and hence globally (see \cite{fmq}, Theorem 4).


\section{Proofs of Theorem \ref{main1} and Theorem  \ref{main2}: sharpness}\label{sharpness}

Fix a point $x\in M$. For simplicity, for the rest of this section we will use the notation
 \be d=d(x,\cdot),\qquad V(r)=V_{x}(r)=\mu(B(x,r)),\qquad  \Lambda(r)=\Lambda_x(r),\ee
and we will assume that $(M,g)$ is not isometric to $\R^n$, in which case $\Lambda$ is strictly decreasing to 0.

Later we will need the following estimate

\be \Lambda(r)\le \wt\L(r)=\min_{\delta>0}\big(\delta+\delta^{-2n}\L(\delta^{2n+1}r)\big).\label{LL}\ee
This holds since $\L$ is decreasing, $0<\L\le 1-\sm$ and
\be \delta+\delta^{-2n}\L(\delta^{2n+1}r)\ge\bc \L(r) & {\text { if } } 0<\delta\le1 \cr 1 & {\text{ if }} \delta>1.\ec\ee

The Euclidean  volume growth condition can be written as 
\be\frac{\sm B_n}{V(r)}\le\frac 1 {r^n}= \frac{\sm B_n}{V(r)}+\frac{B_n\;\L(r)}{ V(r)}\le \frac{\sm B_n}{V(r)}+B_n\frac{\L\big(B_n^{-\frac 1 n}V^{\frac1 n}(r)\big)}{ V(r)}.\label {est3}\ee

 We will also repeatedly use that for $a,r>0$
 \be{a^n}\sm V(r)\le V(ar)\le \frac{a^n}\sm V(r),\label{volestimate}\ee
and
\be a^{-\frac1{2n+1}}\min\{a,1\}\Lt(r)\le  \Lt(ar)\le a^{-\frac1{2n+1}}\max\{a,1\} \Lt(r).\label{Lest}\ee

Assume first $n\ge 3$. Following the notation used in \cite{cm1}, \cite{cm2},   we  consider the function
\be b(y)=\bigg(\frac\sm{c_2} \; G_2(x,y)\bigg)^{\frac1{2-n}},\qquad y\in M\setminus \{x\}\label {b}.\ee
which should be thought of  as a  good, smooth replacement of the distance function, for large distances.

In particular, under our hypothesis the main estimates \eqref{gsharp1} and \eqref{gradgsharp1} for $\a=2$ can be written in terms of $b$ as follows:
\be |b-d|\le \Cl{b} \;d\; \wt\L(d),\qquad d=d(x,y)\ge 1\label{b1}\ee
\be\bigg(\fint_{A_{r,(1+\eta)r}}\big|\nabla b- \nabla d\big|^2 \bigg)^{\frac1 2}\!\!\le \C(1+\eta)^{n-1}\Big(1+\frac{1}{\sqrt\eta}\Big)
 \wt\L^{\frac12}(r),\qquad r\ge r_1\label{gradb1}\ee
for a suitable $r_1>1$ and for any fixed $\eta>0$.

Note that, as pointed out in \cite{cm2} and \cite{c2}, by the Laplacian comparison and the maximum principle, one can show that $G_2(x,y)\ge c_2 d^{2-n}(x,y)$, which implies $b\le \sm^{\frac1{2-n}}d$. In \cite{c2} Colding proved that we also have $|\nabla b|\le \sm^{\frac1{2-n}}$, with equality at one point if and only if $M$ is the flat $\R^n$.

\begin{rk} \label{cheng} In   \cite{cm2}, Remark 2.11, it was remarked that the critical sets  $\{b=r\}\cap\{\nabla b=0\}$ are closed sets of codimension 2, by results of Cheng in \cite{che}, Thm. 2.2, which also show that  $\{b=r\}\setminus \{\nabla b=0\}$ is an $(n-1)-$dim. smooth submanifold of $M$. Hence each set $\{b=r\}$ has zero measure in $M$.
	\end{rk}

Inequality \eqref{b1} follows from \eqref{gsharp1}, with $\alpha=2$, combined with the rough global bound for $G_2$,  given in  \eqref{globalG}, which can be written in terms of $b$ as 
\be \Cl{b1} \;d\le b\le \Cl{b2} \;d,\qquad d=d(x,y)>0.\label{b2}\ee

To prove \eqref{gradb1}, we have
  \be\nabla b=-\frac\sm {(n-2)c_2}b^{n-1}\nabla G_2\label{b6}\ee
 and  from \eqref{globalG}, \eqref{globalgradg} we get $|\nabla b|\le \C.$ 
Also, using \eqref{b6} estimate  \eqref{gradgsharp2} can be rewritten as
\be\bigg(\fint_{A_{r,(1+\eta)r}}\big|b^{1-n}\nabla b- d^{1-n}\nabla d\big|^2 \bigg)^{\frac1 2}\!\!\le \C\Big(1+\frac{1}{\sqrt\eta}\Big)
 \wt\L^{\frac12}(r) r^{1-n}.
 \ee
 From \eqref{b1} we find $\Cl{a1},\Cl{a2}>0$ and $r_1>1$ such that 
 \be d^{-n}\big(1-\Cr{a1}\wt\L(d)\big)\le b^{-n}\le d^{-n}\big(1+\Cr{a2}\wt\L(d)\big),\qquad d=d(x,y)\ge r_1\label{b3}\ee
and using $|b^{1-n}-d^{1-n}|\le \C d^{1-n}\wt\L(d)$ for $d\ge r_1$ large enough, we obtain \eqref{gradb1}.

Additionally, from \eqref{b3} and \eqref{est3}, \eqref{LL}, and the fact that $\Lambda(r)\le 1-\sm$, we find $a_1,a_2>0$ such that 

\be\frac{B_n\sm}{V(d)}\Big(1-a_1 \wt\L\big(B_n^{-\frac 1n}V^{\frac1n}(d)\big)\Big)\le b^{-n}\le \frac{B_n\sm}{V(d)}\Big(1+a_2 \wt\L\big(B_n^{-\frac 1n}V^{\frac1n}(d)\big)\Big),\qquad d\ge r_1.\label{bn}
\ee 
\medskip
 For any $r>0$ let $\rho$ be uniquely related to $r$ as follows
 \be V(r)=V(\rho)\;\exp\Big({\Lt^{-\frac12}\big(V^{\frac1n}(\rho)\big)}\Big)
 \label{rrho}
\ee
(note that $V(r)$ is strictly increasing in $r$). Clearly $r>\rho$ and $\rho\to+\infty$ if $r\to+\infty$, and from \eqref{volestimate}, \eqref{Lest} we also have
\be \Lt^{\frac12}(\rho)\log\frac r\rho\le \Cl{x1} \Lt^{\frac12}\big(V^{\frac1n}(\rho)\big)\log\frac {V(r)}{V(\rho)}+\Cl{x2}= \Cr{x1}+\Cr{x2}.\label{Lest1}\ee
 
For   $t\in\R, r>0,\a>0$ define, with $\rho$ satisfying \eqref{rrho}, 
 \be h_{\a,r}(t)=\bc \rho^{-\a}-r^{-\a}& {\text {if } } \;t\le\rho\\ t^{-\a}-r^{-\a} & {\text {if }} \;\rho<t\le r\\ 0& {\text {if }} \; t>r \ec\qquad \qquad h_{0,r}(t)=\bc   \log\dfrac{\raisebox{-.5ex}{$r$}}{\raisebox{.7ex}{$\rho$}}& {\text {if } } \;t\le\rho\\  \log\dfrac{\raisebox{-.45ex}{$r$}}{\raisebox{.25ex}{$t$}} & {\text {if }} \;\rho<t\le r\\ 0& {\text {if }} \; t>r \ec\ee

Consider the family of functions
\be u_{\a,r}(y)=T_\a\big(h_{\a,r}(b)\big)(y)=\int_M G_\a(y,z)h_{\a,r}(b(z))dz\label{T}\ee

for $\alpha$ even, and for $\alpha$ odd
\be \tilde u_{\a,r}=\bc \dfrac{u_{\a-1,r}}{\a-1}& {\text {if  }}\a\ge 3 \\
h_{0,r}(b) & {\text {if  }}\a=1 \ec\ee
 which clearly satisfy 
 \be\bc D^\a u_{\a,r}=h_{\a,r}(b) & {\text{if }} \alpha\; {\text {even }}\\  D^\a\tilde  u_{\a,r}=\dfrac{\nabla h_{\a-1,r}(b)}{\a-1}&{\text{if }} \alpha\; {\text { odd,}}\;\a\ge3\\D^\a\tilde  u_{1,r}=\nabla h_{0,r}(b) &{\text{if }} \a=1.\ec\ee
 
\medskip
 We will prove that for $r$ large enough
 
 \begin{align} &\|D^\a u_{\a,r}\|_{n/\a}^{n/\a}=B_n \sm \log \frac{V(r)}{V(\rho)}+O(1)\label{sh1}\\[.2em]
  &\|D^\a\tilde u_{\a,r}\|_{n/\a}^{n/\a}=B_n \sm \log \frac{V(r)}{V(\rho)}+O(1)\label{sh2}\\[.2em]
 &u_{\a,r}(y)\ge c_\a B_n \log  \frac{V(r)}{V(\rho)}-C,\qquad y\in B(x,\rho)\label{sh3}\\[.2em]
 &\tilde u_{\a,r}(y)\ge \tilde c_\a B_n \log  \frac{V(r)}{V(\rho)}-C,\qquad y\in B(x,\rho)\label{sh4}\\[.2em]
 &\|u_{\a,r}\|_{n/\a}^{n/\a}\le CV(r),\qquad \a<\frac n2,\label {sh5}\\
  &\|\tilde u_{\a,r}\|_{n/\a}^{n/\a}\le CV(r),\qquad \a<\frac n2+1,\label {sh6}
 \end{align}
 where ``$O(1)$" denotes a function of $r$ (depending on the fixed $x$) which is bounded as $r\to\infty$. 
  
From \eqref{b2} we find $a_3,a_4>0$, $\;a_3<1<a_4$, such that 
 \be B(x,a_3r)\subseteq \{b\le r\}\subseteq B(x,a_4r),\qquad r>0.\label{b4}\ee
 
 \begin{rk} Using \eqref{b1} and \eqref{b2} it's possible to refine the above inclusions \eqref{b4}, for large $r$, and to  approximate $\{b\le r\}$ with balls of radii $r(1\pm a\wt\L(r))$.\end{rk}
 
 \medskip
 

\ni\underline{{Proof of \eqref{sh1}, \eqref{sh2}}}. Let  $\alpha$ be even. If $r$ so large so that $a_3\rho\ge r_1$, using \eqref{b4} we have
 \be\ba&\int_M|D^\a u_{\a,r}|^{\na}=\int_{M} (h_{\a,r}(b))^{\na}=\int_{b\le \rho}(\rho^{-\a}-r^{-\a})^{\na}+\int_{\rho<b\le r} (b^{-\a}-r^{-\a})^{\na}\cr&\le V(a_4)+\int_{a_3\rho<d\le a_4 r} b^{-n} \le V(a_4)+\int_{a_3\rho<d\le a_4 r}  \frac{B_n\sm}{V(d)}\Big(1+a_2 \wt\L\big(B_n^{-\frac 1n}V^{\frac1n}(d)\big)\Big)\\
 &\le V(a_4)+ B_n\sm\int_{V(a_3\rho)}^{V(a_4r)} \frac{1}{t}\Big(1+a_2\wt\L\big(B_n^{-\frac 1n}t^{\frac1n}\big)\Big)dt=V(a_4)+ B_n\sm\log\frac{V(a_4r)}{V(a_3\rho)}\cr&+B_n\sm a_2 \Lt\big(B_n^{-\frac1n}V^{\frac1n}(a_3\rho)\big)\log\frac{V(a_4r)}{V(a_3\rho)}\le 
 B_n\sm\log\frac{V(r)}{V(\rho)}\cr&+
\C \Lt\big(V^{\frac1n}(\rho)\big)\log\frac{V(r)}{V(\rho)}+\C\le B_n\sm\log\frac{V(r)}{V(\rho)}+\C
 \ea\ee
 where we used  \eqref{Lest1}.
 
In the other direction, for $r$ large enough, and arguing as above,
 \be\ba&\int_{M} (h_{\a,r}(b))^{\na}\ge\int_{\rho<b\le r} b^{-n}\bigg(1-\C\Big(\frac b r\Big)^\a\bigg)\\&\ge\int_{a_4\rho<d\le a_3 r} \frac{B_n\sm}{V(d)}\Big(1-a_1 \wt\L\big(B_n^{-\frac 1n}V^{\frac1n}(d)\big)\Big)-\C r^{-\a}\int_{d\le a_3r} d^{\a-n} \cr&\ge B_n\sm\log\frac{V(r)}{V(\rho)}-\Cr{rr}-\C r^{-\a}\int_{d\le a_3r} V^{\frac{\a-n}n}(d)\cr&=
B_n\sm\log\frac{V(r)}{V(\rho)}-\Cl{rr}-\C r^{-\a}\int_0^{V(a_3r)} t^{\frac{\a-n}n}dt
 \label{b9}
 \ea\ee
 which gives \eqref{sh1}.
 
 Now let $\alpha$ be odd. If $\alpha=1$ we get
\be\int_M|D^\a\tilde u_{1,r}|^\na=\int_M |\nabla h_{0,r}(b)|^n=\int_{\rho<b\le r} b^{-n}=B_n\sm\log\frac{V(r)}{V(\rho)}+O(1).
\ee 
If $\alpha\ge3$, 
 \be\ba\int_M|D^\a\tilde u_{\a,r}|^\na=\int_M \Big|\frac{\nabla h_{\a-1,r}(b)}{\a-1}\Big|^\na&=\int_{\rho<b\le r} |\nabla b^{-\a+1}|^{\na}=\int_{\rho<b\le r} b^{-n}|\nabla b|^{\na}
 \\ &=\int_{\rho<b\le r} b^{-n}+\int_{\rho<b\le r} b^{-n}\big(|\nabla b|^{\na}-|\nabla d|^{\na}\big).
 \label{b5}\ea
 \ee 
 
To estimate the last integral, let $J_r\in\N$ be such that $\rho 2^{J_r-1}\le r< \rho 2^{J_r}$, and using \eqref{gradb1}  (with $\eta=1$ and  $r=\rho 2^j$) we get

 \be\ba&\int_{\rho<b\le r} b^{-n}\big||\nabla b|^{\na}-|\nabla d|^{\na}\big|\le \Cl{bc}\int_{\rho<b\le r} b^{-n}|\nabla b-\nabla d| \le \Cr{bc}\sum_{j=0}^{J_r-1} \int_{\rho 2^{j}\le d\le \rho 2^{j+1}}
b^{-n}|\nabla b-\nabla d|\\
& \le \Cr{bc}\sum_{j=0}^{J_r-1} \bigg(\int_{\rho  2^{j}\le d\le \rho  2^{j+1}}
 |\nabla b-\nabla d|^2\bigg)^{1/2}\bigg(\int_{\rho  2^{j}\le d\le \rho  2^{j+1}} b^{-2n}\bigg)^{1/2}\\&
\le\Cl{bc1}\sum_{j=0}^{J_r-1} \wt\L^{\frac12}(\rho  2^j)V^{\frac12}(\rho  2^{j+1}) \bigg(\int_{\rho  2^{j}\le d\le \rho  2^{j+1}} V^{-2}(d)\bigg)^{1/2}\le\Cr{bc1} \sum_{j=0}^{J_r-1} \wt\L^{\frac12}(\rho  2^j)\bigg(\frac{V(\rho 2^{j+1})}{V(\rho 2^j)}\bigg)^{1/2} \\&\le \Cl{zz} \sum_{j=0}^{J_r-1} \wt\L^{\frac12}(\rho  2^j)= \Cr{zz}\wt\L^{\frac12}(\rho )+\frac{\Cr{zz}}{\log 2}\sum_{j=1}^{J_r-1}\int_{\rho  2^{j-1}}^{\rho 2^j} \wt\L^{\frac12}(\rho 2^j)\;\frac{dt}t\le \Cl{b9}+\Cl{b10}\sum_{j=1}^{J_r-1}\int_{\rho  2^{j-1}}^{\rho 2^j} \wt\L^{\frac12}(t)\;\frac{dt}t\cr 
&\le   \Cr{b9}+\Cr{b10}\int_{\rho }^r \wt\L^{\frac12}(t)\;\frac{dt}t\le \Cl{b9}+\Cl{b10}\;\Lt^{\frac12}(\rho)\log \frac r\rho\le \C.
\ea\label{b6a}
 \ee
which, combined with \eqref{b5}, gives \eqref{sh2}.
 
 \smallskip
 
 \ni\underline{{Proof of \eqref{sh3}, \eqref{sh4}}}.
 Let $\alpha $ be even. For $r$ large enough, and with $d=d(x,z)$, using \eqref{gsharp1} and \eqref{b3} we obtain
 \be\ba u_{\a,r}(y)&\ge\int_{\rho<b\le r}G_\a(y,z)(b^{-\a}-r^{-\a})dz\ge  
 \int_{\rho<b\le r}G_\a(x,z) b^{-\a}dz\cr&\hskip4em  +\int_{\rho<b\le r}\big(G_\a(y,z)-G_\a(x,z)\big)b^{-\a}dz-\int_{b\le r}G_\a(y,z) r^{-\a}dz\cr&
 \ge\int_{a_4\rho\le d\le a_3r} c_\a\sm^{-1} b^{-n}(1-\C\wt\L(d)\big)dz -\C\int_{ d\le a_3r}  d^{\a-n}(y,z) r^{-\a}dz\cr&\hskip3em -\C\int_{a_4\rho\le d\le a_3r}|G_\a(y,z)-G_\a(x,z)|d^{-\a}(x,z)dz.\label{b7}
 \ea\ee
 
 \smallskip
Since $a_4>1$,   using the Lipschitz estimate \eqref{holderG} and arguing as in \eqref{b9}, using \eqref{bn},  we get that for $y\in B(x,\rho)$

 \be\ba u_{\a,r}(y)&\ge c_\a B_n\log \frac{V(r)}{V(\rho)}-\C-\C\int_{d\ge a_4\rho} d(x,y) d^{-n-1}(x,z)dz\\&\ge c_\a B_n\log \frac{V(r)}{V(\rho)}-\C. \ea\ee

 If $\alpha $ odd  and $\alpha\ge3$ then \eqref{sh4} follows  since $\tilde c_\a=c_\a/(\a-1)$, whereas if $\alpha=1$ then for $y\in B(x,\rho)\subseteq \{b\le \rho/a_3\}$, then 
 \be \tilde u_{0,r}(y)=\bc \log (r/\rho) &{\text{ if }} b\le \rho\\ \log  (r/b) & { \text{ if }} \rho<b\le r\ec \ge\frac1 n \log \frac{V(r)}{V(\rho)}-\Cl{0}=\tilde c_1 B_n \log \frac{V(r)}{V(\rho)}-\Cr{0}.\ee

\medskip  \ni\underline{{Proof of \eqref{sh5}, \eqref{sh6}}}.
Let $\alpha$ be even.  If $d(x,y)\ge 2a_4r$ then
 \be\ba&|u_{\a,r}(y)|=\int_{b\le r}G_\a(y,z)h_{\a,r}(b)(z)dz\le\int_{b\le r} G_\a(y,z) b^{-\a}(z)dz\\&\le \C\int_{d(x,z)\le a_4 r}d^{\a-n}(y,z) d^{-\a}(x,z)dz\le\C d(x,y)^{\a-n}\int_{d(x,z)\le a_4 r} d^{-\a}(x,z)\le\C d^{\a-n}(x,y) r^{n-\a}.\ea\ee
If $\a<n/2$ then $d^{\a-n}(x,\cdot)\in L^{\na}(B(x,R)^c)$ for any $R>0$, and 
 \be\int_{d(x,y)\ge 2a_4r}|u_{\a,r}(y)|^{\na}dy\le\C r^{(n-a)\na}\int_{d(x,y)\ge 2a_4r} d(x,y)^{(\a-n)\na}dx\le \C r^n\le \C V(r).\ee
 
 Now, due to \eqref{globalG}, the operator $T_\a$ in \eqref{T} is bounded from $L^p(M)$ to $L^q(M)$ with $1<p<\frac n\a$ and $q^{-1}=p^{-1}-\frac{\a}{n}$ (this is the Hardy-Littlewood -Sobolev inequality for generalized, Riesz-like  potentials, see \cite{a2}).
 So if we pick any $p$ such that $\frac n{2\a}<p<\na$ then $q>\na$, and  for $R=2a_4r$
 \be\ba&\bigg( \int_{d(x,y)\le R}|u_{\a,r}(y)|^{\na}dy\bigg)^{\an}\le \C V(R)^{\an-\frac1q}\bigg( \int_{d(x,y)\le R}|T_\a\big(h_{\a,r}(b)\big)|^{q}dy\bigg)^{\frac1q}\\&\le \C V(R)^{\an-\frac1q}\bigg(\int_{d(x,y)\le R}\big|h_{\a,r}(b)\big|^p\bigg)^{\frac1p}\le \C V(R)^{\an-\frac1q}R^{\frac n p-\a}\le 
 \C V(R)^\an\le \C V(r)^\an,
 \ea
 \ee
 which concludes the proof of \eqref{sh5} for $\alpha$ even, $\alpha<n/2$, and therefore also for \eqref{sh6} for $\alpha$ odd,  $3\le\alpha<\frac n2+1.$ When $\alpha=1$ one checks easily that \eqref{sh6} also holds.
 
 \medskip
Let us now prove the sharpness statement. Suppose that $\alpha$ is even, $\a<\frac n 2$,  $n\ge 3$, and consider $v_{\a,r}=u_{\a,r}/\|D^\a u_{\a,r}\|_{n/\a}$. From \eqref{sh1}, \eqref{sh3}, \eqref{sh5}, we get, for $r$ large enough,
 \be|v_{\a,r}(y)|^\nna\ge c_\a^{\nna}B_n \sm^{-\frac\a{n-\a}}\log\frac {V(r)}{V(\rho)}-C,\qquad y\in B(x,\rho)\label{g1}\ee
 \be\|D^\a v_{\a,r}\|_{n/\a}=1,\qquad \|v_{\a,r}\|_{n/\a}^{n/\a}\le C \frac{V(r)}{\log \frac{V(r)}{V(\rho)}},\label{g2}\ee
 and recall that from \eqref{rrho}
 \be \log\frac{V(r)}{V(\rho)}=\frac1{\Lt^{\frac12}\big(V^{\frac1n}(\rho)\big)}\to+\infty,\qquad r\to\infty.\ee
 Hence, for $\theta>1$ and $r$ large enough, recalling that $\gamma_{n,\a}^{-1}=c_\a^\nna B_n$, and from estimate (3.2) in  \cite{q},
 \be\ba\frac{1}{\|v_{\a,r}\|_{n/\a}^{n/\a}}&\int_M\frac{\exp_{\llceil\frac{n-\a}{\a}-1\rrceil}\left(\theta\sm^{\frac\a{n-\a}}\gamma_{n,\a}|v_{\a,r}|^\nna\right)}{1+|v_{\a,r}|^\nna}dy\cr&\ge \frac{1}{\|v_{\a,r}\|_{n/\a}^{n/\a}}\int_{|v_{\a,r}|\ge1}\frac{\exp\left(\theta\sm^{\frac\a{n-\a}}\gamma_{n,\a}|v_{\a,r}|^\nna\right)}{1+|v_{\a,r}|^\nna}dy-C\cr&\ge \frac{1}{\|v_{\a,r}\|_{n/\a}^{n/\a}}\int_{B(x,\rho)}\frac{\exp\left(\theta\sm^{\frac\a{n-\a}}\gamma_{n,\a}|v_{\a,r}|^\nna\right)}{1+|v_{\a,r}|^\nna}dy-C\cr&\ge C\frac{V(\rho)}{V(r)}\log \frac{V(r)}{V(\rho)}\cdot \frac{\exp\big(\theta\log\frac{ V(r)}{V(\rho)}-C\big)}{\log \frac{V(r)}{V(\rho)}}-C\cr&=C \bigg(\frac{V(r)}{V(\rho)}\bigg)^{\theta-1}-C\to+\infty,\qquad r\to+\infty
 \ea \label{sharp}\ee
 which gives the sharpness of the exponential constant in \eqref{MS}, for $\alpha$ even.  Regarding the sharpness of the power 
 in the denominator of \eqref{MS}, using the same estimates as in \eqref{sharp} we get for $\theta<1$
 
  \be\ba\frac{1}{\|v_{\a,r}\|_{n/\a}^{n/\a}}&\int_M\frac{\exp_{\llceil\frac{n-\a}{\a}-1\rrceil}\left(\sm^{\frac\a{n-\a}}\gamma_{n,\a}|v_{\a,r}|^\nna\right)}{1+|v_{\a,r}|^{\theta\nna}}dy\cr&\ge \frac{1}{\|v_{\a,r}\|_{n/\a}^{n/\a}}\int_{|v_{\a,r}|\ge1}\frac{\exp\left(\sm^{\frac\a{n-\a}}\gamma_{n,\a}|v_{\a,r}|^\nna\right)}{1+|v_{\a,r}|^{\theta\nna}}dy-C\cr&\ge \frac{1}{\|v_{\a,r}\|_{n/\a}^{n/\a}}\int_{B(x,\rho)}\frac{\exp\left(\sm^{\frac\a{n-\a}}\gamma_{n,\a}|v_{\a,r}|^\nna\right)}{1+|v_{\a,r}|^{\theta\nna}}dy-C\cr&\ge C\frac{V(\rho)}{V(r)}\log \frac{V(r)}{V(\rho)}\cdot \frac{\exp\big(\log\frac{ V(r)}{V(\rho)}-C\big)}{\big(\log \frac{V(r)}{V(\rho)}\big)^\theta}-C\cr&=C \bigg(\log \frac{V(r)}{V(\rho)}\bigg)^{1-\theta}-C\to+\infty,\qquad r\to+\infty.
 \ea \label{sharpden}\ee
 
 The sharpness for $\alpha$ odd, $\a<\frac n 2+1$,  follows in the same exact manner, using $v_{\a,r}=\tilde u_{\a,r}/\|D^\a\tilde u_{\a,r}\|_{n/\a}$, and \eqref{sh2}, \eqref{sh4}, \eqref{sh6}.
 
 \smallskip
 It remains to consider the cases $\a=1, n=2,$ and $\a=2, n=3,4.$. For $\a=1$ we note that  when $n\ge3$ we  could have just as well considered the functions 
 $\tilde u_{1,r}(y)=h_{0,r}\big(d(x,y)\big)$, which are defined also for  $n=2$, since they are Lipschitz, compactly supported, hence in $W^{1,n}(M)$, and satisfy the same estimates as in \eqref{sh2}, \eqref{sh4}, \eqref{sh6}, and therefore \eqref{sharp}, for the corresponding $v_{1,r}$.

 When $\alpha=2$ the proof of sharpness for arbitrary $n\ge3$ can be effected as follows. Let $h_r$ be smooth on $[0,\infty)$, decreasing, such that $|h'_r|\le C,$ $|h_r''|\le C$ with 
\be h_r(t)=\bc \log\dfrac{\raisebox{-.65ex}{$r$}}{\raisebox{.55ex}{$\rho$}}  & \text { if } 0<t<\rho\\\log\dfrac{\raisebox{-.55ex}{$r$}}{\raisebox{.25ex}{$t$}} & {\text {if }} \;\rho+1<t\le r-1\\ 0& {\text {if }} \; t>r, \ec\label{fam}\ee
 and  $h_r(t)\le \log (r/t)$ on $[r-1,r]$. We then let 
  $u_r=h_r(b)$, which belongs to $W_0^{2,n/2}(M)$. One checks that $\Delta b=(n-1)b^{-1}|\nabla b|^2$, and therefore
\be \Delta u_r=\Delta h_r(b)=\big(h_r''(b)+h_r'(b)(n-1)b^{-1}\big)|\nabla b|^2\label{deltahr}\ee
and in particular
\be \Delta u_r=(2-n)b^{-2}|\nabla b|^2,\qquad \rho+1<b< r-1.\ee
Using the same estimates  as in \eqref{b5}, \eqref{b6a} we get for $r $ large enough
\be\ba\int_M|\Delta u_r|^{\frac n2}&=(n-2)^{\frac n2}\int_{\rho<b\le r} b^{-n}|\nabla b|^{n}+O(1)\\&=
(n-2)^{\frac n2}\int_{\rho<b\le r} b^{-n}+(n-2)^{\frac n2}\int_{\rho<b\le r} b^{-n}\big(|\nabla b|^{n}-|\nabla d|^{n}\big)+O(1)\\&=(n-2)^{\frac n2}B_n\sm \log \frac{V(r)}{V(\rho)}+O(1)\big),
\ea
\ee
where in the first identity above we used $|\nabla b|\le C$ and 
\be\int_{r-1<b\le r} b^{-n}\le C,\ee
and the same estimate with $\rho$ rather than $r$. 

Hence,
\be\|\Delta u_r\|_{n/2}=\frac{\sm^{2/n}}{n c_2B_n^{\frac{n-2}n}}\bigg(\log\frac {V(r)}{V(\rho)}\bigg)^{\frac2n}+O(1)\big)\ee and furthermore, if $d(x,y)\le \rho$ then $b\le \rho/a_3<r-1$, for $r$ large enough, so for $y\in B(x,\rho)$
\be u_r(y)\ge\bc \log (r/\rho) &{\text{ if }} b< \rho\\ \log  (r/b) & { \text{ if }} \rho\le b\le \rho/a_3\ec\ge \log \frac r\rho-C\ge \tilde c_1 B_n\log\frac{V(r)}{V(\rho)}-C .\ee
 \be \|u_r\|_{n/2}^{n/2}\le C \log r+\int_{2<b\le r} \Big(\log\frac r b\Big)^{n/2}\le C V(r).\ee
 Then, the function $v_r=u_r/\|\Delta u_r\|_{n/2}$ satisfies the same estimates as $v_{2,r}$ in \eqref{g1}, \eqref{g2}, \eqref{sharp}, which gives the  sharpness statements for \eqref{MS}, when $\alpha=2$ and for any $n\ge 3$.
 
 \medskip
 \begin{rk}\label{b} Except for the case $\alpha=2$, $n=3,4$, the proof of sharpness goes through using $d(x,y)$ instead of $b(y)$. We decided to use $b$ throughout, since the estimates needed for $\alpha=2$, using the family $\{u_r\}$,  are basically the same as those used for $\alpha>2$ using the family $\{u_{\a,r}\}$  (e.g.  \eqref{b5}, \eqref{b6a}). We were not able to use the family $\{u_r\}$ above, to settle the sharpness for  $\alpha>2$.
 \end{rk}

 \begin{rk}\label{rrho1}  The fine calibration between $\rho$ and $r$ given in \eqref{rrho} is only needed to prove sharpness of the power in the denominator of \eqref{MS}. Indeed, the proof of the sharpness of the exponential constant in \eqref{MS} goes through using the same families of functions as above, but with $\rho=1$.
 
 \end{rk}

  \medskip
 To prove  the sharpness statement for \eqref{MSR} and \eqref{MT}  consider the function $h_r$ satisfying \eqref{fam}, and such that $|h_r^{(j)}|\le C$, $j=1,...,\a$. For a fixed $r_0$ s.t. $0<r_0<\inj(M)$ let
 \be w_\e(y)=h_{1/\e}\Big(\frac{d(x,y)}{\e r_0}\Big)\ee
 which is nothing but the manifold version of the standard Adams-Moser family, i.e. a smoothing of the function $\log(r_0/d(x,y))$, as defined on $2\e r_0<d(x,y)\le r_0-\e r_0$, for $\e$ small enough. With calculations as in \cite{f} one finds that for $\alpha$ even and $\e$ small enough

\begin{align} &\|D^\a w_\e\|_{n/\a}=\frac1{n c_\a B_n^{\frac{n-\a}n}}   \Big(\log \frac1{V(\e)}\Big)^{\frac\a n}+O(1)\label{sh1a}\\[.1em]
   &|w_\e(y)|\ge \frac1n\log \frac1{V(\e)}-C,\qquad y\in B(x,\e r_0)\label{sh2a}\\[.2em]
  &\|w_\e\|_{n/\a}\le C\label {sh3a}
 \end{align}
 and the same estimates hold for $\alpha$ odd, with $\ct_\a$ in place of $c_\a$. The sharpness of the exponential constants in \eqref{MSR}, \eqref{MT}, and of the power in the denominator of \eqref{MSR},   follows easily from the above estimates, considering the family $v_\e=w_\e/(\kappa\|w_\e\|_{n/\a}^{n/\a}+\|D^\a w_\e\|_{n/\a}^{n/\a})^{\a/n}$.
 
 \qed
 

\appendix

\section{Proof of \eqref{reg1}-\eqref{reg4}}

\vskip1em

For a fixed $t$ let  $f(r)=E(r,t)$ so that 
$$f'(r)=-\frac{r}{2t} E(r,t),\qquad f''(r)=\Big(-\frac{1}{2t}+\frac{r^2}{4t^2}\Big) E(r,t).$$
Note that 
$$|f'(r)|=\frac{r}{2(4\pi)^{n/2}}\, t^{-n/2-1}\exp\Big(\!\!-\frac{r^2}{4t}\Big)\le \Cl{A1} r^{-n-1}.$$
Hence using the MVT, for each $y$ there is $\tilde d_j$ with $|\tilde d_j-d|\le |d_j-d|<1/j$ s.t. for $j>\frac1 2\dist(x,K)$
$$|E(d_j,t)-E(d,t)|=|d_j-d| f'(\tilde d_j) \le \frac{\C}{j} \Big(d-\frac1 j\Big)^{-n-1}\le \frac \C j \dist(x,K)^{-n-1}.$$

We have 
$$\nabla E(d_j,t)=-\nabla d_j \frac{d_j}{2t} E(d_j,t)\to-\nabla d \frac{d}{2t} E(d,t)=\nabla E(d,t),$$
and
$$|\nabla E(d_j,t)|\le  \frac {d_j}{t} E(d_j,t)\le \C \dist (x,K)^{-n-1}.$$

Using that $\Delta (f\circ g) =f'(g)\Delta g+f''(g)|\nabla g|^2$, and using  iv), we get
\be\ba \Delta E(d_j,t)&=\bigg(-\frac{d_j}{2t}\Delta d_j-\frac{1}{2t}+\frac{d_j^2}{4t^2}\bigg) E(d_j,t)\cr&\ge
\bigg(-\frac{n-1}{2t}\frac{d_j}{d}-\frac{d_j}{2jt}-\frac{1}{2t}+\frac{d_j^2}{4t^2}\bigg) E(d_j,t)\cr&\ge
\bigg(-\frac{n-1}{2t}\bigg|\frac{d_j}{d}-1\bigg|-\frac{d_j}{2jt}-\frac{n}{2t}+\frac{d_j^2}{4t^2}\bigg) E(d_j,t)\cr&
\ge-\bigg(\frac{n-1}{2j d}+\frac{d_j}{2j}\bigg) t^{-n/2-1}\exp\Big(\!\!-\frac{d_j^2}{4t}\Big)-\frac{n}{2}t^{-n/2-1}\cr&
\ge-\frac{\C}{j}\big(\dist(x,K)^{-n-3}+\dist(x,K)^{-n-1}\big)-\frac{n}{2}t^{-n/2-1}.
\ea\ee

\begin{table}[h]

\setlength{\tabcolsep}{24pt} 
\begin{tabular}{@{}p{0.5\linewidth}p{0.5\linewidth}@{}}
\textbf{Luigi Fontana} & \textbf{Carlo Morpurgo} \\
Dipartimento di Matematica ed Applicazioni & Department of Mathematics \\
Universit\'a di Milano-Bicocca & University of Missouri \\
Milan, 20125 & Columbia, Missouri 65211 \\
Italy & USA \\
\texttt{luigi.fontana@unimib.it} & \texttt{morpurgoc@umsystem.edu} \\\\
\textbf{Liuyu Qin} & \\
Department of Mathematics and Statistics & \\
Hunan University of Finance and Economics & \\
Changsha, Hunan & \\
China & \\
\texttt{Liuyu\_Qin@outlook.com} & \\
\end{tabular}
\end{table}
\end{document}